\documentclass[11pt]{article}
\usepackage{graphicx}
\usepackage{amsmath}
\usepackage{amsfonts}
\usepackage{epsfig}
\usepackage{setspace}
\usepackage{diagbox}
\usepackage{authblk}
\usepackage{caption}
\usepackage{subcaption}

\newcommand{\be}{\begin{equation}}
\newcommand{\ee}{\end{equation}}
\newcommand{\beqn}{\begin{eqnarray}}
\newcommand{\eeqn}{\end{eqnarray}}
\newcommand{\beqns}{\begin{eqnarray*}}
\newcommand{\eeqns}{\end{eqnarray*}}

\newcommand{\card}{\mbox{card}}
\newcommand{\Var}{\mbox{Var}}

\newcommand{\EE}{\ensuremath{{\mathbb E}}}
\newcommand{\II}{\ensuremath{{\mathbb I}}}

\newcommand{\fr}[1]{(\ref{#1})}

\newcommand{\Te}{\Theta}

\newtheorem{lemma}{Lemma}
\newtheorem{theorem}{Theorem}

\newtheorem{remark}{Remark}

\begin{document}

\title{\Large{\bf Adaptive estimation for the nonparametric bivariate additive model in random design with long-memory dependent errors}}

\author[1]{Rida Benhaddou\thanks{benhaddo@ohio.edu}}
\author[2]{Qing Liu }
\affil[1]{Department of Mathematics, Ohio University, Athens, OH 45701}
\affil[2]{Department of Mathematics and Statistics, Wake Forest University, Winston Salem, NC 27109}

\date{}

\doublespacing
\maketitle
\begin{abstract}
We investigate the nonparametric bivariate additive regression estimation in the random design and long-memory errors  and construct adaptive thresholding estimators based on wavelet series. The proposed approach achieves asymptotically near-optimal convergence rates when the unknown function and its univariate additive components  belong to Besov space. We consider the problem under two noise structures; (1) homoskedastic Gaussian long memory errors and (2) heteroskedastic Gaussian long memory errors. In the homoskedastic long-memory error case, the estimator is completely adaptive with respect to the long-memory parameter. In the heteroskedastic long-memory case, the estimator may not be adaptive with respect to the long-memory parameter unless the heteroskedasticity is of polynomial form. In either case, the convergence rates depend on the long-memory parameter only when long-memory is strong enough, otherwise, the rates are identical to those under i.i.d. errors. The proposed approach is extended to the general $r$-dimensional additive case, with $r>2$, and the corresponding convergence rates are free from the curse of dimensionality. \\


{\bf AMS (2000) Subject Classification: 62G05, 62G20, 62G08 }
 \end{abstract} 

\section{Introduction}

Consider the nonparametric additive regression model  
\be
y_i=y(t_i, x_i) =\beta_o+ f(t_i)+g(x_i) +\sigma(t_i, x_i) \varepsilon_i, \ \ \ i=1, 2, \cdots, N, \label{conveq}
\ee
where  $t_i$ and $x_i$ are i.i.d. random variables with known compactly supported marginal probability density functions $h_1$ and $h_2$, respectively, and $t_i$ and $x_i$ are independent for any  $i\in \{1, 2, \cdots, N\}$. $\sigma(t,x)$ is (1) a known positive variance constant, that is $\sigma(t,x)=\sigma$ and (2) $\sigma(t,x)$ is a known non-constant deterministic bivariate function. The bivariate additive regression function $U(t, x)=\beta_o+f(t)+g(x)$ is unknown. The univariate functions $f$ and $g$ are real-valued, integrate to zero and are ${\bf L}^2\left([0,1]\right)$, and $\beta_o $ is the overall mean parameter which is a real, nonzero constant. $\{ \varepsilon_i\}_{i\geq1}$ is a stationary Gaussian sequence that is independent of $t_i$ and $x_i$ for any  $i\in \{1, 2, \cdots, N\}$. Furthermore, $\{\varepsilon_i\}$ suffers from long-memory which is described by the linear structure
\be
\varepsilon_i=\sum^{\infty}_{m=0}a_m\eta_{i-m},\ with\ a_o=1,
\ee
where $\{\eta_i\}_{i \in Z}$ is an i.i.d. Gaussian sequence and
\beqns
\lim_{m\rightarrow \infty} a_m m^{(\alpha+1)/2}=1,\ for\ \alpha \in (0, 1),
\eeqns
and
\be \label{longm}
\Var\left(\sum^N_{i=1}\varepsilon_i\right)\sim \pi_{\alpha}N^{2-\alpha},
\ee
where $\pi_{\alpha}$ is a finite and positive constant, and $\alpha$ is the long-memory parameter. The goal is to estimate the additive regression function $U(t, x)=\beta_o+ f(t)+g(x)$ and its components $\beta_o$, $f(t)$ and $g(x)$ based on the data points $(t_1, x_1, y_1)$, $(t_2, x_2, y_2)$, $\cdots$, $(t_N, x_N, y_N)$.   A model of this sort is referred to as the nonparametric additive regression or just additive model.

Additive models are advantageous in their ability to balance the strengths of the fully nonparametric and parametric estimation procedures. In particular, additive models tend to have estimators with lower variances than fully nonparametric models, and lower bias than their parametric counterparts. This problem has a great deal of applications, in particular, {the analysis of responses that are non-linearly associated with several predictors. In such models, the response depends linearly on unknown univariate functions of individual predictor variables and the goal is to make inference about the individual effects and the global effect of these predictors through these functions.

This problem, under various settings, has been studied considerably by the means of a number of nonparametric methods, including kernel smoothing, local polynomials and wavelets among others, and the list of articles includes, in chronological order, Stone~(1985), Buja, Hastie and Tibshirani~(1989), Linton and Nielsen (1995), Opsomer and Ruppert~(1997) and Amato and Antoniadis~(2001), to name a few. Recently, Chesneau, Fadili and Maillot~(2015) investigated a tensor product wavelet solution to the problem under random design, weakly dependent predictors and i.i.d. errors. The main difference between the current work and theirs is that our predictors are independent from one another, and the errors suffer from long-memory under both the homoskedasticity and heteroskedasticity. Grez and Vidakovic~(2018) considered the problem in random design and  i.i.d. errors setting, and proposed another wavelet approach to solve it. The key difference between our work and their work is that their joint design density belongs to some generalized Holder class, it is unknown and needs to be estimated from the data, and their predictors ${\bf x}= \left(x_1, x_2, \cdots, x_p\right)$ may not be independent. Most recently, Amato et al.~(2022) proposed a wavelet-based estimation and model selection for the nonparametric additive model under i.i.d. Gaussian and sub-Gaussian errors. The major difference between the present paper and their work is that they suggest an estimation procedure for the additive components and a selection procedure of the sparse additive components when the number of additive components is allowed to grow with the sample size, so their work is technically high-dimensional.    All of the articles above assume that the error terms are white noise processes or { i.i.d.} noise. More precisely, to the best of our knowledge, no article has studied long-memory in the context of additive models under random design, so we intend to fill in the gap. 

Long-memory has been investigated quite considerably in many nonparametric estimation problems, including regression and deconvolution and the list includes Wang~(1996, 1997), Csorgo and Mielniczuk~(1999), Beran and Feng~(2001), Yang~(2001), Comte, Dedecker and Taupin~(2008),  Kulik and Raimondo~(2009), Kulik and Wichelhaus~(2011), Wishart~(2013), Benhaddou, Kulik, Pensky and Sapatinas~(2014), Benhaddou~(2016, 2018a, 2018b, 2021), and Benhaddou and Liu~(2019). 

The treatment of nonparametric wavelet estimation under irregular or random design has been considered in the literature and one of the first to address nonparametric wavelet estimation under irregularly spaced data points was Cai and Brown~(1998, 1999). This was followed by several articles, but we will only list a few and they are in chronological order Csorgo and Mielniczuk~(1999), Pensky and Vidakovic~(2001), Yang~(2001), Chesneau~(2007), Kulik and Raimondo~(2009), Antoniadis, Pensky and Sapatinas~(2014), Chesneau, Fadili and Maillot~(2015), Grez and Vidakovic~(2018), and Benhaddou~(2020). 

The objective of the paper is to solve the nonparametric bivariate additive regression model with both homoskedastic and heteroskedastic long-memory Gaussian noise via wavelet hard-thresholding when the predictors are random, independent from one another and follow known probability density functions $h_1$ and $h_2$. Motivated to some extent by the work of Kulik and Raimondo~(2009), we derive lower bounds for the $L^2$-risk when the univariate components $f$ and $g$ of the additive regression function belong to some Besov class, and with the help of some results for martingales, we construct adaptive wavelet-thresholding estimators  for $f$ and $g$. In addition, we show that the proposed estimator attains asymptotically quasi-optimal convergence rates in the minimax sense. Furthermore, we demonstrate that long-memory has a detrimental effect on the convergence rates only when it is strong enough  relative to the classical rates obtained under the i.i.d. setting, in which case, the convergence rates depend only on the long-memory parameter $\alpha$. A limited simulations study illustrates the finite sample performance of the proposed estimator under several heteroskedasticity scenarios and long-memory levels. The proposed approach is extended to the general $r$-dimensional additive case, with $r>2$ when the number of additive components is not very large, and therefore the problem is not high-dimensional. Under such setup, it turns out that the corresponding convergence rates are free from the curse of dimensionality. The present work is also motivated by  Grez and Vidakovic~(2018). It turns out that if we fix one of the univariate additive functions $f$ (or $g$) or set it equal to zero, our convergence rates are identical to those that correspond to the dense case in Kulik and Raimondo~(2009) for their study of warped wavelet nonparametric univariate regression under random design and similar noise structure. So, our work can be viewed as an extension of theirs to the nonparametric additive model with the exception that the wavelet coefficients are estimated differently.   Moreover, even in the special case of  long-memory parameter $\alpha=1$, our convergence rates are quite different from  Grez and Vidakovic~(2018), in their treatment of the problem when the functions under consideration belong to certain Sobolev class and errors are i.i.d. This is mainly due to the fact that their joint design density belongs to some generalized Holder class, it is not known, and is to be estimated from the data. Nevertheless, when their joint design density is uniform over the torus $[0, 1]^r$, their convergence rates depend only on the smallest smoothness parameter from amongst the univariate functions, and therefore their rates do not suffer from the curse of dimensionality. Thus, this behavior is consistent with ours. Finally, with long-memory parameter $\alpha=1$, our convergence rates are similar to those in Chesneau et al.~(2015). 

 \section{Estimation Algorithm}

Let $\psi$ and $\eta$ be wavelet bases with compact support and $s_{1o}$ and $s_{2o}$ vanishing moments, respectively, and let $m_{1o}$ and $m_{2o}$ be the lowest resolution levels. Let the scaling functions for $\psi$ and $\eta$ be denoted by $\psi^{(s)}$ and $\eta^{(s)}$, respectively. Since $f, g \in {\bf L}^2\left([0, 1]\right)$, they can be extended as wavelet series  
\be
f(t)=\sum^{2^{m_{1o}}-1}_{k_1=0}\theta^{(1)}_{m_{1o}, k_1}\psi^{(s)}_{m_{1o}, k_1}+\sum^{\infty}_{j_1=m_{1o}}\sum^{2^{j_1}-1}_{k_1=0}\beta^{(1)}_{j_1,k_1}\psi_{j_1, k_1}(t),
\ee
\be
g(x)=\sum^{2^{m_{2o}}-1}_{k_2=0}\theta^{(2)}_{m_{2o}, k_1}\eta^{(s)}_{m_{2o}, k_2}+\sum^{\infty}_{j_2=m_{2o}}\sum^{2^{j_2}-1}_{k_2=0}\beta^{(2)}_{j_2,k_2}\eta_{j_2, k_2}(x),
\ee 
where $\theta^{(1)}_{j_1, k_1}=\int^1_0f(t)\psi^{(s)}_{j_1, k_1}(t)dt$, $\theta^{(2)}_{j_2, k_2}=\int^1_0g(x)\eta^{(s)}_{j_2, k_2}(x)dx$, $\beta^{(1)}_{j_1, k_1}=\int^1_0f(t)\psi_{j_1, k_1}(t)dt$ and $\beta^{(2)}_{j_2, k_2}=\int^1_0g(x)\eta_{j_2, k_2}(x)dx$. Here, set the lowest resolution levels equal to zero, that is, $m_{io}=0$, $i=1, 2$. \\
Then, inspired by  Grez and Vidakovic~(2018), allow the unbiased estimators for the wavelet and scaling coefficients
\be \label{betest1}
\widehat{\beta^{(1)}_{j_1, k_1}}=\frac{1}{N}\sum^N_{i=1}\frac{y(t_i, x_i)}{h_1(t_i)h_2(x_i)}\psi_{j_1k_1}(t_i),
\ee
\be   \label{betest2}
\widehat{\beta^{(2)}_{j_2, k_2}}=\frac{1}{N}\sum^N_{i=1}\frac{y(t_i, x_i)}{h_1(t_i)h_2(x_i)}\eta_{j_2k_2}(x_i),
\ee
\be  \label{thetest1}
\widehat{\theta^{(1)}_{j_1, k_1}}=\frac{1}{N}\sum^N_{i=1}\frac{y(t_i, x_i)}{h_1(t_i)h_2(x_i)}\left(\psi^{(s)}_{j_1k_1}(t_i)-2^{-j_1/2}\right),
\ee
\be   \label{thetest2}
\widehat{\theta^{(2)}_{j_2, k_2}}=\frac{1}{N}\sum^N_{i=1}\frac{y(t_i, x_i)}{h_1(t_i)h_2(x_i)}\left(\eta^{(s)}_{j_2k_2}(x_i)-2^{-j_2/2}\right).
\ee
Therefore, define the hard thresholding estimators for $f$ and $g$,
\be\label{fest}
\widehat{f}_N(t)=\widehat{\theta^{(1)}_{0,0}}\psi^{(s)}_{0,0}(t)+\sum^{J_1-1}_{j_1=m_{1o}}\sum^{2^{j_1}-1}_{k_1=0}\widehat{\beta^{(1)}_{j_1,k_1}}\II\left(|\widehat{\beta^{(1)}_{j_1k_1}}|> \lambda^{(1)}(j_1)\right)\psi_{j_1, k_1}(t),
\ee
\be \label{gest}
\widehat{g}_N(x)=\widehat{\theta^{(2)}_{0, 0}}\eta^{(s)}_{0,0}(x)+\sum^{J_2-1}_{j_2=m_{2o}}\sum^{2^{j_2}-1}_{k_2=0}\widehat{\beta^{(2)}_{j_2,k_2}}\II\left(|\widehat{\beta^{(2)}_{j_2k_2}}|> \lambda^{(2)}(j_2)\right)\eta_{j_2, k_2}(x),
\ee
where, the quantities $J_i$, $m_{io}$, and $\lambda_i(j_i)$, $i=1,2$, will be determined under the two different setups in the proceeding sections. In addition, an unbiased estimator for $\beta_o$ is given by
\be \label{conset}
\widehat{\beta}_o=\frac{1}{N} \sum^N_{i=1} \frac{y(t_i, x_i)}{h_1(t_i)h_2(x_i)}.
\ee
Thus, the estimator for the bivariate function $U(t, x)$ is the sum of its components' own estimators \fr{conset}, \fr{fest} and \fr{gest}, as follows
\be \label{uest}
\widehat{U}_N(t, x)=\widehat{\beta}_o+\widehat{f}_N(t)+\widehat{g}_N(x).
\ee
 Next is the list of conditions that will be utilized in the derivation of the theoretical results. \\
{\bf Assumption A.1.}\label{A1} $f, g \in {\bf{L^{2}}}\left[0, { 1} \right]$ are both bounded above, that is, there exists positive constant $M_1$, $M_2<\infty$ such that  $f(t) \leq M_1$ and $g(x)\leq M_2$, for all  $t, x \in \left[0, { 1} \right]$.\\
{\bf Assumption A.2.} \label{A2} The probability density functions $h_1$ and $h_2$ are uniformly bounded on $\left[0, {1} \right)$, that is, there exist positive constants $m_{11}$, $m_{12}$, $m_{21}$ and $m_{22}$, with $0< m_{i1}\leq m_{i2} < \infty$, $i=1,2$, such that $m_{11} \leq h_1(t) \leq m_{12}$ and $m_{21} \leq h_2(x) \leq m_{22}$. 
 \begin{remark}	  
{ Assumption {\bf A.2.} is valid for instance when  $h_1$ and $h_2$ are the uniform distributions. }
  \end{remark}
{\bf Assumption A.3.}  \label{A3}To ensure unique identification of the additive components, let the functions $f$ and $g$ be such that $\int^1_0f(t)dt=0$ and $\int^1_0g(x)dx=0$. \\
  {\bf Assumption A.4.} \label{A4}Denote for $i=1,2$,
  \beqn  \label{eq10}
 s^{*}_i&=&s_i+1/2 - 1/p.
 \eeqn
 The univariate functions $f(t)$ and $g(x)$ belong to a Besov space. In particular, if $s_{io} \geq s_i$, $i=1,2$, their wavelet coefficients $ \beta^{(i)}_{j_ik_i}$ satisfy
\be  \label{eq11}
 B^{s_1, s_2}_{p, q}(A)=\left \{ U \in {\bf L}^2([0, 1]): \left( \sum_{j_i} 2^{(j_is_i^{*})q}\left (\sum_{k_i}| \beta^{(i)}_{j_i, k_i}|^{p}\right)^{q/{p}}\right )^{1/q} \leq A_i, i=1, 2\right \}.
\ee
 We are in the position to fill in the details of the estimator and find the minimax lower bound for the quadratic risk and compare it to asymptotic upper bound for the mean squared error of our estimator. We define the minimax $L^2$-risk over a set $\Theta$ as 
$$R(\Theta)=\inf_{\tilde{f}}\sup_{f\in \Theta}\EE\|\tilde{f}-f\|^2,$$
where the infimum is taken over all possible estimators $\tilde{f}$ of $f$. 
 \section{Asymptotic minimax and adaptivity: $\sigma(t, x)\equiv\sigma$ constant case}
 According to Yang~(2001), the minimax risk of estimating a regression function under random design and long-memory errors is the sum of two terms; the minimax rate of the same class but under $i.i.d.$ errors, and the rate of estimating the mean of the regression function under correlated errors. So, we will use this result in the derivation of the lower bounds. 
 \begin{theorem}\label{th:lowerbds} Let $\min\{s_1, s_2\} \geq \max\{\frac{1}{p}, \frac{1}{2} \}$ with $1 \leq p,q \leq \infty$, and $A_i > 0$, $i=1,2$. Let  Assumptions {\bf{A.1-A.4}} hold, and let condition \fr{longm} hold. Then, for $i=1, 2$, as $N\rightarrow \infty$,
 \be \label{lowerbds}
 R({\bf{B^{s_1, s_2}}}(A))\geq C\left\{ \begin{array}{ll} 
A_2^2 \left[\frac{1}{A_2^2N} \right]^{\frac{2s_2}{2s_2+1}}\vee \frac{1}{N^{\alpha}} , & \mbox{if}\ \  {s_1}> {s_2},\\
A_1^2 \left[\frac{1}{A_1^2N}\right]^{\frac{2s_1}{2s_1 +1}}\vee \frac{1}{N^{\alpha}}   , & \mbox{if}\ \   {s_1}\leq {s_2}.
   \end{array} \right.
   \ee
 \end{theorem}
 \begin{remark}	  
{ Notice that in line with Yang~(2001) argument, if $h_1$ and $h_2$ are uniform, then $\EE(f(t))=\int^1_0f(t)h_1(t)dt=\int^1_0f(t)dt=0$, and $\EE(g(x))=\int^1_0g(x)h_2(x)dx=\int^1_0g(x)dx=0$, and therefore, the minimax risk for estimating $f$ and $g$ will not involve the long-memory parameter $\alpha$. However, the minimax risk for estimating $U(t, x)=\beta_o+f(t)+g(x)$ will, since $\EE(U(t, x))=\beta_o\neq0$. }
  \end{remark}
 \begin{lemma} \label{lem:Var}
Let assumptions ${\bf{A.1}}$, ${\bf{A.2}}$ and ${\bf{A.3}}$ hold and let $\widehat{\beta^{(i)}_{j_ik_i}}$ and $\widehat{\theta^{(i)}_{j_ik_i}}$ be defined in \fr{betest1}-\fr{thetest2}. Then, under condition \fr{longm}, as $N \rightarrow \infty$, one has 
\be \label{var-bet}
\EE\|\widehat{\beta^{(i)}_{j_ik_i}}-\beta^{(i)}_{j_ik_i}\|^2 \asymp \frac{1}{N}.
\ee
\be \label{var-thet}
\EE\|\widehat{\theta^{(i)}_{j_ik_i}}-\theta^{(i)}_{j_ik_i}\|^2 \asymp \frac{1}{N}.
\ee
\be \label{var2-bet}
\EE\|\widehat{\beta^{(i)}_{j_ik_i}}-\beta^{(i)}_{j_ik_i}\|^4 \asymp \frac{1}{N^2}.
\ee
In addition, for the estimator $\widehat{\beta}_o$ in \fr{conset}, under condition \fr{longm} as $N \rightarrow \infty$, one has 
\be \label{var-beto}
\EE|\widehat{\beta}_{o}-\beta_{o}|^2 \asymp \frac{\sigma^2}{N^{\alpha}}.
\ee
\end{lemma}
Based on {\bf Lemma 1}, choose the thresholds $\lambda^{(i)}(j_i)$, $i=1, 2$, such that 
\be \label{thresh1}
\lambda^{(i)}(j_i)=\gamma\frac{{\ln(N)}}{\sqrt{N}}.
\ee
 In addition, the highest resolution levels $J_1$ and $J_2$ should be chosen such that
\be \label{J}
2^{J_1}= 2^{J_2}\sim \frac{N}{\ln(N)}.
\ee 
  \begin{remark}	  
Notice that under homoskedastic long-memory errors, the estimators of the additive components $f$ and $g$ are completely adaptive with respect to the parameter of the long memory $\alpha$. This behavior is consistent with the literature for univariate regression case under similar error structure (e.g., Kulik and Raimondo~(2009)). However, the convergence rates of estimating the bivariate additive regression as a whole may depend on such quantity. 
  \end{remark}
 \begin{lemma} \label{lem:Lardev}
Let $\widehat{\beta^{(i)}_{j_ik_i}}$, $i=1, 2$, be defined in \fr{betest1} and \fr{betest2}, respectively, and let assumptions ${\bf{A.1}}$-${\bf{A.3}}$ and condition \fr{longm} hold. Then, as $N \rightarrow \infty$, one has 
\be  \label{Largdev}
\Pr \left(| \widehat{\beta^{(i)}_{j_ik_i}}-\beta^{(i)}_{j_ik_i} |> 1/2\lambda({j_i})\right)\asymp \left[\frac{1}{N}\right]^{2\tau},
\ee
where $\tau$ is a positive parameter that is large enough. 
\end{lemma}
\begin{theorem} \label{th:upperbds-2}
Let  $\widehat{U}_{N}(t, x)$ be the estimator defined in \fr{uest} with its components $\widehat{\beta}_o$, $\widehat{f}_N(t)$ and $\widehat{g}_N(x)$ given in \fr{conset}, \fr{fest} and \fr{gest}, respectively. Suppose assumptions ${\bf{A.1}}$-${\bf{A.4}}$ hold. Then, under condition \fr{longm}, if $\tau$ is large enough, as $N \rightarrow \infty$, one has
 \be \label{upperbds-2}
 \sup_{U\in B^{s_1, s_2}(A)} \EE\|\widehat{U}_N-U\|^2\leq C\left\{ \begin{array}{ll} 
 A_2^2 \left[\frac{\ln(N)}{A_2^2N} \right]^{\frac{2s_2}{2s_2+1}}\vee \frac{1}{N^{\alpha}} , & \mbox{if}\  {s_1}\geq {s_2},\\
  A_1^2 \left[\frac{\ln(N)}{A_1^2N} \right]^{\frac{2s_1}{2s_1  +1}}\vee \frac{1}{N^{\alpha}} , & \mbox{if} \   {s_1}< {s_2}.\\
  \end{array} \right.
\ee
\end{theorem}
 \begin{remark}	
{\bf Theorem \ref{th:lowerbds}} and {\bf Theorem \ref{th:upperbds-2}} imply that, for the $L^2$-risk, estimator \fr{uest} with its components $\widehat{\beta}_o$, $\widehat{f}_N(t)$ and $\widehat{g}_N(x)$ given in \fr{conset}, \fr{fest} and \fr{gest}, respectively, is adaptive and asymptotically near-optimal within a logarithmic factor of $N$ over all  Besov balls $B^{s_1, s_2}(A)$. 
  \end{remark} 
  \begin{remark}	  
One key difference between the present work and that of Grez and Vidakovic~(2018) is that the joint design density in their work is unknown and needs to be estimated from the data, and the predictors ${\bf x}= \left(x_1, x_2, \cdots, x_p\right)$ are not necessarily independent from one another. In the special case of  $\alpha=1$, our convergence rates are quite different from theirs, in their treatment of the problem when the functions under consideration belong to certain Sobolev space, the joint design density belongs to some generalized Holder class of functions, and errors are i.i.d. This is mainly due to the fact that their joint design density is unknown and to be estimated from the data. 
  \end{remark}
   \begin{remark}	
Our work can be viewed as an extension of Kulik and Raimondo~(2009) to the nonparametric additive model with the exception that the wavelet coefficients are estimated differently. In addition, if we fix the univariate function $g$ (or $f$) or set it equal to zero, our rates will be identical to theirs, for their dense case in their consideration of nonparametric univariate regression estimation under random design and the same error structure based on warped wavelets. Furthermore, the sparse case does not contribute to our convergence rates at all. This behavior is not consistent with their work, but it is consistent with Chesneau et al.~(2015).
\end{remark}
 \begin{remark}	
Our work is comparable to that of Chesneau et al.~(2015). In the special case of  $\alpha=1$, our rates are identical to theirs, in their consideration of additive regression estimation under random design and moderately dependent predictors based on tensor product wavelet basis. 
\end{remark}
  \section{Asymptotic minimax and adaptivity: $\sigma(t, x)$ non-constant case}
\noindent
{\bf Assumption A.5.}  $\sigma(t, x)$ is bounded below on $[0, 1]^2$, that is, there exist constant $\sigma_0>0$ such that for all $t$ and $x$ in $[0, 1]$, $\sigma(t, x) >\sigma_o$.
\begin{theorem}\label{th:lowerbds-lm} Let  Assumptions {\bf{A.1-A.5}} hold. Then, as $N\rightarrow \infty$, one has
 \be \label{lowerbds-lm}
 R({\bf{B^{s_1, s_2}}}(A))\geq C\left\{ \begin{array}{ll} 
A_2^2 \left[\frac{1}{A_2^2N} \right]^{\frac{2s_2}{2s_2+1}}\vee \frac{1}{N^{\alpha}} , & \mbox{if}\ \  {s_1}> {s_2},\\
A_1^2 \left[\frac{1}{A_1^2N}\right]^{\frac{2s_1}{2s_1 +1}}\vee \frac{1}{N^{\alpha}}   , & \mbox{if}\ \   {s_1}\leq {s_2}.
   \end{array} \right.
   \ee
 \end{theorem}
  \begin{lemma} \label{lem:Var-lm}
Let conditions ${\bf{A.1}}$-${\bf{A.3}}$ hold and let $\widehat{\beta^{(i)}_{j_ik_i}}$ and $\widehat{\theta^{(i)}_{j_ik_i}}$ be defined in \fr{betest1}-\fr{thetest2}. Then, provided that $\int^1_0\int^1_0\psi_{j_1k_1}(t)\sigma(t, x)dtdx\neq 0$ and $\int^1_0\int^1_0\eta_{j_2k_2}(x)\sigma(t, x)dtdx\neq 0$, as $N \rightarrow \infty$, one has 
\be \label{var-bias}
\EE\|\widehat{\beta^{(i)}_{j_ik_i}}-\beta^{(i)}_{j_ik_i}\|^2 \asymp \frac{1}{N}+\frac{1}{N^{\alpha}}.
\ee
\be \label{var-thet}
\EE\|\widehat{\theta^{(i)}_{j_ik_i}}-\theta^{(i)}_{j_ik_i}\|^2 \asymp \frac{1}{N}+ \frac{1}{N^{\alpha}}.
\ee
\be \label{var2-bet}
\EE\|\widehat{\beta^{(i)}_{j_ik_i}}-\beta^{(i)}_{j_ik_i}\|^4 \asymp \frac{1}{N^2}.
\ee
In addition, for the estimator $\widehat{\beta}_o$ in \fr{conset}, under condition \fr{longm} as $N \rightarrow \infty$, one has 
\be \label{var-beto}
\EE|\widehat{\beta}_{o}-\beta_{o}|^2 \asymp \frac{\sigma^2}{N^{\alpha}}.
\ee
\end{lemma}
Based on {\bf Lemma 3}, choose the thresholds $\lambda^{(1)}(j_1)$ and $\lambda^{(2)}(j_2)$ such that 
\be \label{thresh21}
\lambda^{(1)}(j_1)=\left\{ \begin{array}{ll}  
\gamma\frac{{\ln(N)}}{\sqrt{N}}, \ & \mbox{if}\ \int^1_0\int^1_0\psi_{j_1k_1}(t)\sigma(t, x)dtdx=0, \\
\gamma_1\sqrt{\frac{\ln{(N)}}{N^{\alpha}}}\ & \mbox{if}\ otherwise,
  \end{array} \right.
  \ee
  and
\be \label{thresh22}
\lambda^{(2)}(j_2)=\left\{ \begin{array}{ll}  
\gamma'\frac{{\ln(N)}}{\sqrt{N}}, \ & \mbox{if}\ \int^1_0\int^1_0\eta_{j_2k_2}(x)\sigma(t, x)dtdx=0. \\
\gamma'_1\sqrt{\frac{\ln{(N)}}{N^{\alpha}}}\ & \mbox{if}\ otherwise,
  \end{array} \right.
  \ee
 In addition, the highest resolution levels $J_1$ and $J_2$ should be chosen such that
\be \label{J2}
2^{J_1}= 2^{J_2}\sim \frac{N}{\ln(N)}.
\ee 
  \begin{remark}	  
As you can see from \fr{thresh21}, \fr{thresh22} and \fr{J2}, under heteroskedastic long-memory errors, the estimators for the additive components $f$ and $g$ may not be adaptive with respect to the parameter of the long memory $\alpha$. Nevertheless, if the heteroskedasticity structure $\sigma(t, x)$ is of polynomial form it suffices to use wavelet functions with enough number of vanishing moments to achieve full adaptivity with respect to $\alpha$. This is consistent with Kulik and Raimondo~(2009). 
  \end{remark}
 \begin{lemma} \label{lem:Lardev-lm}
Let $\widehat{\beta^{(i)}_{j_ik_i}}$, $i=1, 2$, be defined in \fr{betest1} and \fr{betest2}, respectively, and let assumptions ${\bf{A.1}}$-${\bf{A.3}}$ and condition \fr{longm} hold. Then, as $N \rightarrow \infty$, one has 
\be  \label{Largdev-lm}
\Pr \left(| \widehat{\beta^{(i)}_{j_ik_i}}-\beta^{(i)}_{j_ik_i} |> 1/2\lambda^{(i)}({j_i})\right)\asymp \left[\frac{1}{N}\right]^{2\tau},
\ee
where $\tau$ is a positive parameter that is large enough. 
\end{lemma}
\begin{theorem} \label{th:upperbds-lm}
Let  $\widehat{U}_{N}(t, x)$ be the estimator defined in \fr{uest} with its components $\widehat{\beta}_o$, $\widehat{f}_N(t)$ and $\widehat{g}_N(x)$ given in \fr{conset}, \fr{fest} and \fr{gest}, respectively. Suppose assumptions ${\bf{A.1}}$-${\bf{A.5}}$ hold. Then, under condition \fr{longm}, if $\tau$ is large enough, as $N \rightarrow \infty$, one has
 \be \label{upperbds-lm}
 \sup_{U\in B^{s_1, s_2}(A)} \EE\|\widehat{U}_N-U\|^2\leq C\left\{ \begin{array}{ll} 
 A_2^2 \left[\frac{\ln(N)}{A_2^2N} \right]^{\frac{2s_2}{2s_2+1}}\vee \frac{1}{N^{\alpha}} , & \mbox{if}\  {s_1}\geq {s_2},\\
  A_1^2 \left[\frac{\ln(N)}{A_1^2N} \right]^{\frac{2s_1}{2s_1  +1}}\vee \frac{1}{N^{\alpha}} , & \mbox{if} \   {s_1}< {s_2},\\
  \end{array} \right.
 \ee
\end{theorem}

\section{Simulation study}
To assess the performance of our estimator defined in \fr{uest} in a finite sample setting and to investigate the effect of the long-memory on the estimator, we carry out a limited simulation study. We evaluate the mean integrated square error (MISE) $\EE \| \hat{y}-y\|^2$ of the  estimator. All simulation is completed with MATLAB (R2022a).

\begin{enumerate}
\item We generate the data using equation \fr{conveq} with test functions $f(t_i) =  8\left(  t_i - \frac{1}{2} \right)^2  - \frac{2}{3}$ and $g(x_i)=- cos(4\pi x_i +1 )$, where  $t_i$ and $x_i$ are i.i.d. random variables from the uniform distributions ($h_1=h_2=1$) defined on $[0, 1]$, $i = 1,2, . . . ,N$, with $N=2^{12}=4096$. For the constant parameter $\beta_0$, we choose $\beta_0 = 0.25* \text{median}(\{f(t_i) +g(x_i) \})$. 

\item We simulate the long-memory error $\{\varepsilon_i\}$ in \fr{conveq} using a one-dimensional fractionally differenced ARIMA (ARFIMA) sequences with the fractional differencing parameter $d \in \{0.1, 0.2, 0.3, 0.4 \}$. Note that the fractional differencing parameter $d$ is obtained from $\alpha$ by $2d={1-\alpha}$. The corresponding choices of $\alpha$ are $\{0.8, 0.6, 0.4, 0.2\}$.
\item The choice of $\sigma$ in \fr{conveq} is determined by the signal-to-noise ratio (SNR), with $\text{SNR}=20\log_{10}\left( \frac{ \| U(t, x)  \|^2   }{ \|\sigma(t,x) \varepsilon \|^2}\right)$. We consider both cases when $\sigma$ is constant and non-constant, and two SNR levels, SNR = 10 dB (medium noise) and SNR =  15 dB (low noise).
\item To compute the hard thresholding estimators \fr{fest} and \fr{gest}, we use the same Daubechies wavelet with 3 vanishing moments, where the thresholds are given in \fr{thresh1}.  And highest resolution levels $J_1$ and $J_2$ are given by \fr{J2}.

\item We compute  the averages of the errors over $n_0=1000$ repeated simulations by
\be
\widehat{\text{MISE}}=\frac{1}{n_0}\sum\limits_{i=1}^{n_0} \left\| \hat{y_i}-y_i \right\|^2.
\ee
\end{enumerate}

\begin{table}[htbp]\centering
\caption{}
MISE averaged over 1000 repeated simulations. \\
Finest resolution levels $J_1=J_2=6$ are used.
\begin{tabular}{|l|cccc|}
\hline
 \multicolumn{5}{|l|}{$\sigma_1(t,x) = $ constant} \\
  \hline
  \diagbox{SNR}{$\alpha$}   & 0.8    & 0.6    & 0.4  &  0.2   \\

\hline
  10dB                 & 0.0247	  &  0.0280	 &  0.0313	 & 0.0321       \\
  15dB                &   0.0104    &  	0.0115    & 	0.0127	    & 0.0130      \\
\hline
\end{tabular}
 
\begin{tabular}{|l|cccc|}
\hline
 \multicolumn{5}{|l|}{$\sigma_2(t,x) =(t+0.4)(x+0.6) $ } \\
  \hline
  \diagbox{SNR}{$\alpha$}   & 0.8    & 0.6    & 0.4  &  0.2   \\

\hline
  10dB                     &0.0256	 &0.0290	& 0.0325 & 	0.0341        \\
  15dB                     & 0.0107 &	0.0119	& 0.0132	& 0.0136  \\
\hline
\end{tabular}

\begin{tabular}{|l|cccc|}
\hline
 \multicolumn{5}{|l|}{$\sigma_3(t,x) =(x+0.6)^2(t+0.4) $ } \\
  \hline
  \diagbox{SNR}{$\alpha$}   & 0.8    & 0.6    & 0.4  &  0.2   \\
\hline
  10dB                     & 0.0253 & 0.0290& 0.0334& 0.0348       \\
  15dB                     & 0.0111& 0.0123& 0.0133 &0.0140  \\
\hline
\end{tabular}

\begin{tabular}{|l|cccc|}
\hline
 \multicolumn{5}{|l|}{$\sigma_4(t,x) =(x+0.6)\sqrt{t+0.4} $ } \\
  \hline
  \diagbox{SNR}{$\alpha$}   & 0.8    & 0.6    & 0.4  &  0.2   \\
\hline
  10dB                     & 0.0248	  &0.0291	  &0.0378	  &0.0826      \\
  15dB                     & 0.0098	  &0.0117	  &0.0195	  &0.0645 \\
\hline
\end{tabular}

\begin{tabular}{|l|cccc|}
\hline
 \multicolumn{5}{|l|}{$\sigma_5(t,x) =(t+0.4)(2+cos(2\pi x)) $ } \\
  \hline
  \diagbox{SNR}{$\alpha$}   & 0.8    & 0.6    & 0.4  &  0.2   \\
\hline
  10dB                     & 0.0230  & 0.0271 &   0.0360    & 0.0790    \\
  15dB                     &  0.0094 &  0.0114 & 0.0181  &0.0612  \\
\hline
\end{tabular}
\end {table}

\begin{figure}
    \caption{}
        \caption*{Functions estimation with $\sigma_5(t,x) =(t+0.4)(2+cos(2\pi x))$. Blue dots represent the actual function $U(t,x)$, red dots represent the estimated function $\widehat{U}(t,x)$.}

  \hspace*{-1.6cm}   \begin{tabular}{llll}
     \begin{subfigure}[b]{0.28\textwidth}
  
     \includegraphics[scale=0.3]{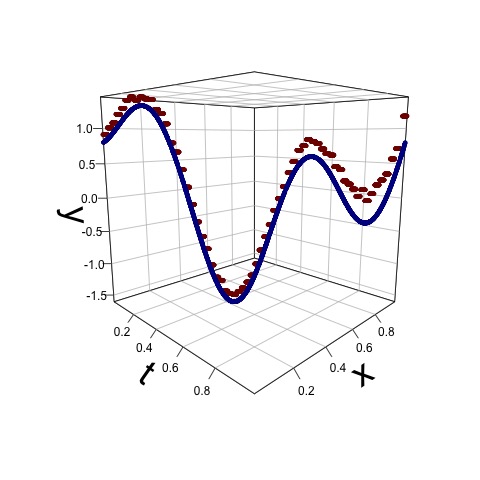}
  \caption*{SNR=10, $\alpha=0.8$}
     \end{subfigure}%
     \hfill
     \begin{subfigure}[b]{0.28\textwidth}
 
        \includegraphics[scale=0.3]{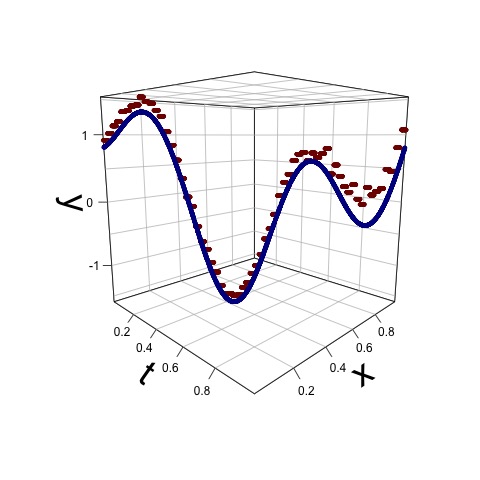}
  \caption*{SNR=10, $\alpha=0.6$}
     \end{subfigure}%
 
     \hfill
     \begin{subfigure}[b]{0.28\textwidth}
 
     \includegraphics[scale=0.3]{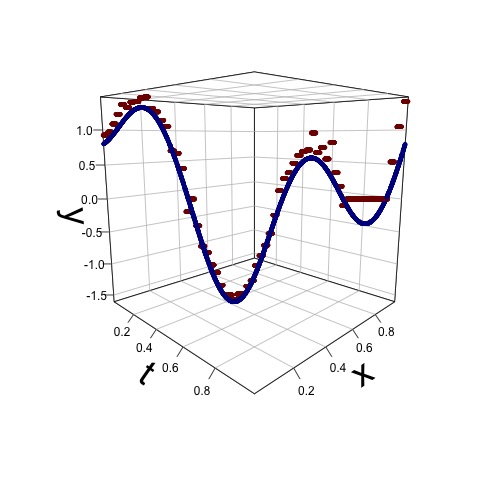}
      \caption*{SNR=10, $\alpha=0.4$}
     \end{subfigure}%
          \hfill
     \begin{subfigure}[b]{0.28\textwidth}
     
          \includegraphics[scale=0.3]{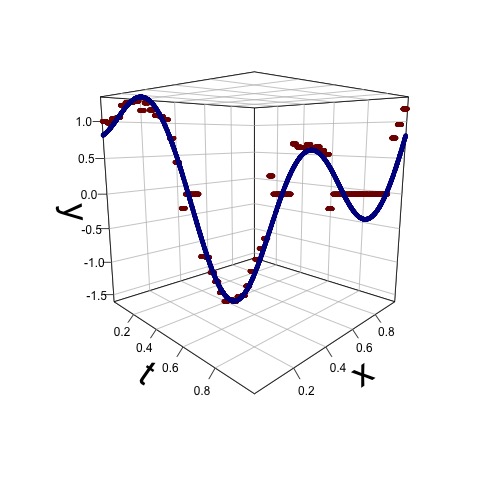}
          \caption*{SNR=10, $\alpha=0.2$}
     \end{subfigure}
     \end{tabular}
      \hspace*{-1.6cm}   \begin{tabular}{llll}
     \begin{subfigure}[b]{0.28\textwidth}
  
     \includegraphics[scale=0.3]{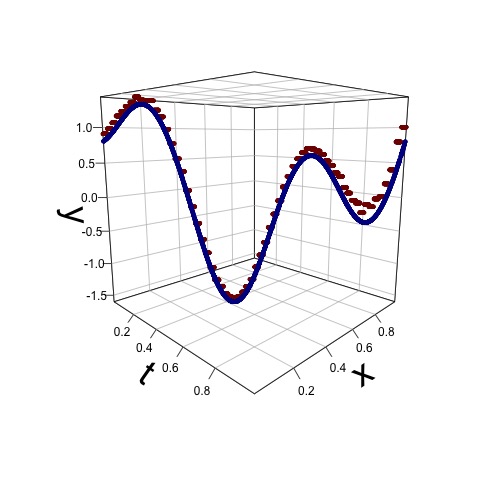}
  \caption*{SNR=15, $\alpha=0.8$}
     \end{subfigure}%
     \hfill
     \begin{subfigure}[b]{0.28\textwidth}
 
        \includegraphics[scale=0.3]{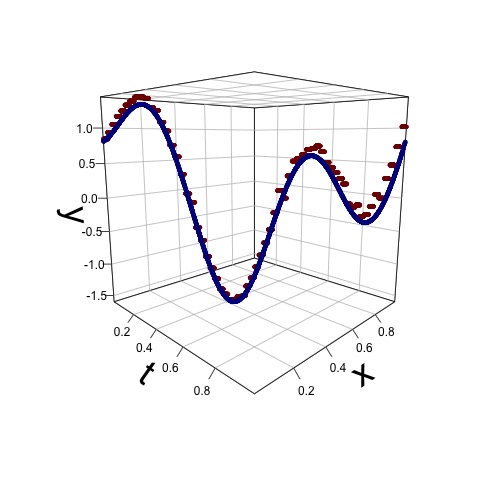}
  \caption*{SNR=15, $\alpha=0.6$}
     \end{subfigure}%
 
     \hfill
     \begin{subfigure}[b]{0.28\textwidth}
 
     \includegraphics[scale=0.3]{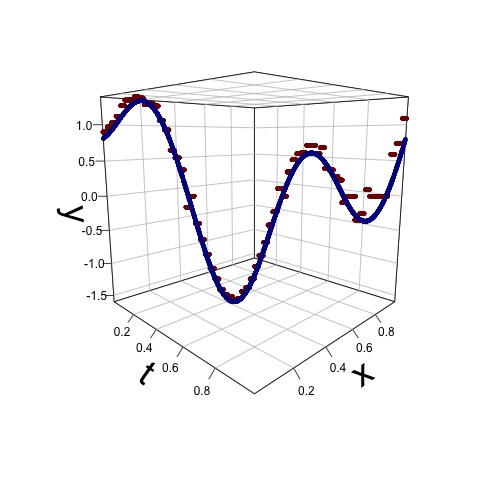}
      \caption*{SNR=15, $\alpha=0.4$}
     \end{subfigure}%
          \hfill
     \begin{subfigure}[b]{0.28\textwidth}
     
          \includegraphics[scale=0.3]{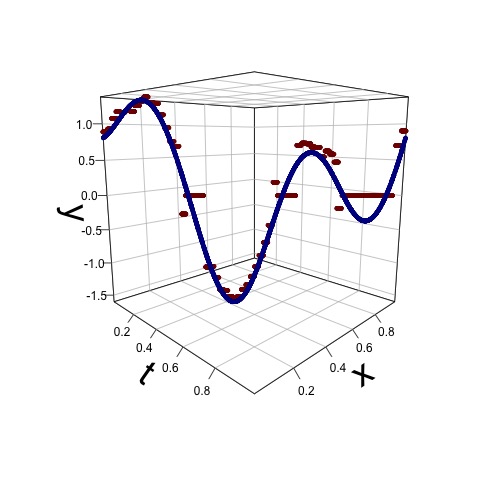}
          \caption*{SNR=15, $\alpha=0.2$}
     \end{subfigure}
     \end{tabular}

        \label{fig1}
\end{figure}

\begin{figure}
    \caption{}
            \caption*{Functions estimation with $\sigma_3(t,x) =(x+0.6)^2(t+0.4) $. Blue dots represent the actual function $U(t,x)$, red dots represent the estimated function $\widehat{U}(t,x)$.}
  \hspace*{-1.6cm}   \begin{tabular}{llll}
     \begin{subfigure}[b]{0.28\textwidth}
  
     \includegraphics[scale=0.3]{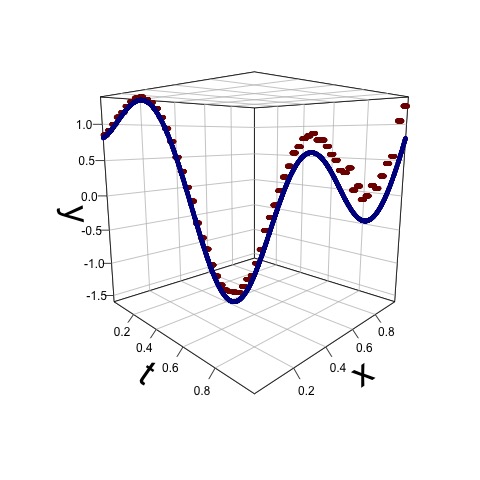}
  \caption*{SNR=10, $\alpha=0.8$}
     \end{subfigure}%
     \hfill
     \begin{subfigure}[b]{0.28\textwidth}
 
        \includegraphics[scale=0.3]{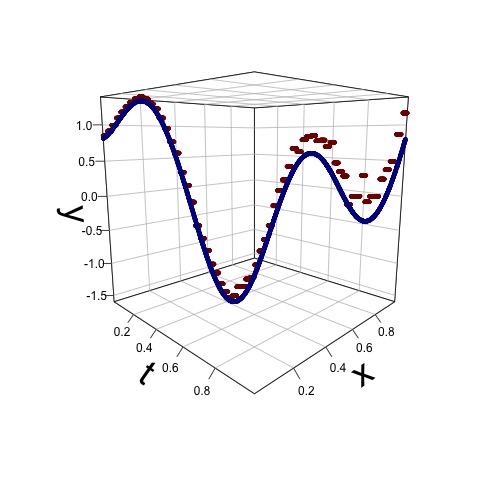}
  \caption*{SNR=10, $\alpha=0.6$}
     \end{subfigure}%
 
     \hfill
     \begin{subfigure}[b]{0.28\textwidth}
 
     \includegraphics[scale=0.3]{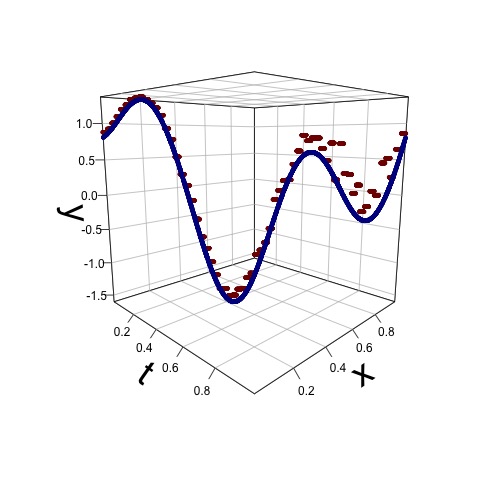}
      \caption*{SNR=10, $\alpha=0.4$}
     \end{subfigure}%
          \hfill
     \begin{subfigure}[b]{0.28\textwidth}
     
          \includegraphics[scale=0.3]{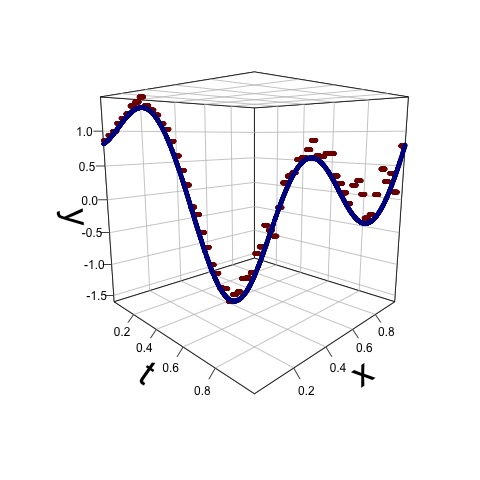}
          \caption*{SNR=10, $\alpha=0.2$}
     \end{subfigure}
     \end{tabular}
      \hspace*{-1.6cm}   \begin{tabular}{llll}
     \begin{subfigure}[b]{0.28\textwidth}
  
     \includegraphics[scale=0.3]{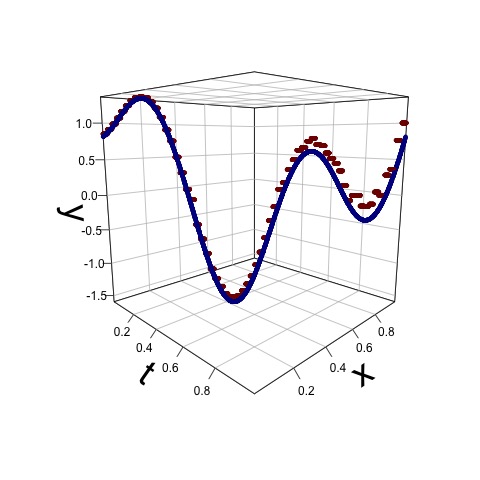}
  \caption*{SNR=15, $\alpha=0.8$}
     \end{subfigure}%
     \hfill
     \begin{subfigure}[b]{0.28\textwidth}
 
        \includegraphics[scale=0.3]{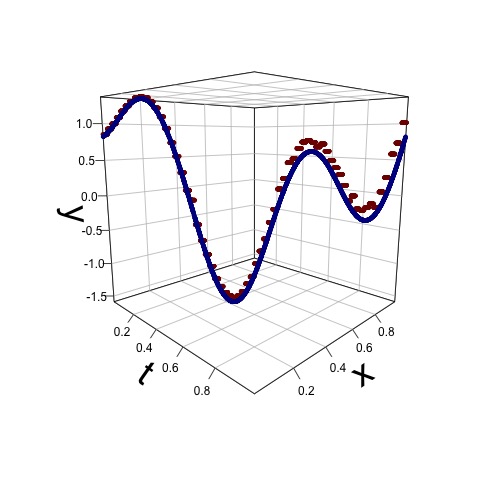}
  \caption*{SNR=15, $\alpha=0.6$}
     \end{subfigure}%
 
     \hfill
     \begin{subfigure}[b]{0.28\textwidth}
 
     \includegraphics[scale=0.3]{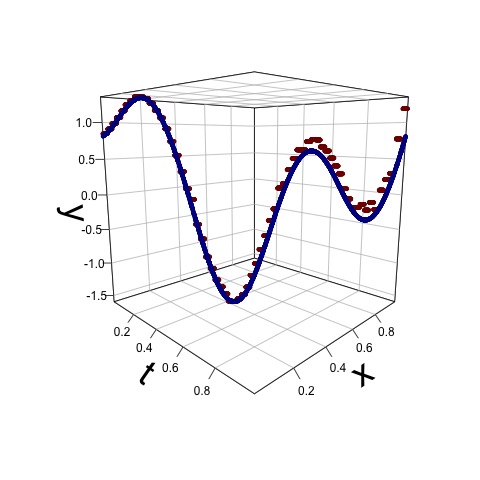}
      \caption*{SNR=15, $\alpha=0.4$}
     \end{subfigure}%
          \hfill
     \begin{subfigure}[b]{0.28\textwidth}
     
          \includegraphics[scale=0.3]{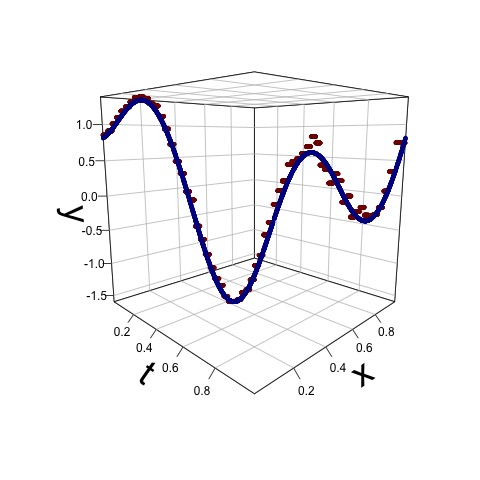}
          \caption*{SNR=15, $\alpha=0.2$}
     \end{subfigure}
     \end{tabular}

        \label{fig2}
\end{figure}

Table 1 reports the averages of the errors over 1000 repeated simulations using a sample size of $N= 2^{12}$, with both homoskedastic and heteroskedastic long-memory errors. It can be seen that, at different signal-to-noise ratios, the estimation algorithm provides a relatively good precision that deteriorates as the long-memory parameter $\alpha$ decreases from 0.8 to 0.2. Figure \fr{fig1} and Figure \fr{fig2} illustrate the function $U(t, x)$ and its estimate $\widehat{U}(t, x)$ and also show a relatively good precision that deteriorates as the LM parameter decreases. It is worth mentioning that estimation errors with both irregular $\sigma(t,x)$ ($\sigma_4$ and $\sigma_5$) are remarkably sensitive to severe long memory ($\alpha= 0.4$ and $\alpha= 0.2$) compared to $\sigma(t,x)$ of polynomial forms ($\sigma_2$ and $\sigma_3$).  The results in Table 1, Figure \fr{fig1} and Figure \fr{fig2} confirm our theoretical findings that as the level of long-memory feature increases ($\alpha$ decreases), the performance deteriorates (MISE increases). 
   \section{Estimation for the general $r$-dimensional additive model, minimax rates and adaptivity}
Consider the general $r$-dimensional additive model, with $r>2$ predictors ${\bf t}= \left(t_1, t_2, \cdots, t_r\right)$, 
\be
y_i=y({\bf t}_i) =\beta_o+ \sum^r_{m=1}f_m(t_{mi}) +\sigma({\bf t}_i) \varepsilon_i, \ \ \ i=1, 2, \cdots, N, \label{conveq2}
\ee
where  $t_{mi}$, $m=1, 2, \cdots, r$, are i.i.d. random variables with known compactly supported probability density functions $h_i$, $m=1, 2, \cdots, r$, and $t_{1i}$, $t_{2i}$, $\cdots$ and $t_{ri}$ are independent for any  $i\in \{1, 2, \cdots, N\}$. $\sigma({\bf t})$ is a known deterministic $r$-dimensional function. The $r$-dimensional additive regression function $U({\bf t})=\beta_o+\sum^r_{m=1}f_m(t_m)$ is unknown. The univariate functions $f_m$, with $m=1, 2, \cdots, r$, are real-valued, integrate to zero and are ${\bf L}^2\left([0,1]\right)$, and $\beta_o $ is a real nonzero constant. The goal is to estimate the additive regression function $U({\bf t})$ and its components $\beta_o$, $f_m$, with $m=1, 2, \cdots, r$, based on the data points $({\bf t}_1, y_1)$, $({\bf t}_2, y_2)$, $\cdots$, $({\bf t}_N, y_N)$. $\{ \varepsilon_i\}_{i\geq1}$ is a stationary Gaussian sequence that is independent of ${\bf t}_i$ for any  $i\in \{1, 2, \cdots, N\}$ and satisfying \fr{longm}.\\
Let $\psi^{(m)}$ be wavelet bases with compact support and $s_{mo}$ vanishing moments, respectively, and let $j_{mo}$ be the lowest resolution levels, with $m=1, 2, \cdots, r$. Let the scaling functions for $\psi^{(m)}$ be denoted by $\psi^{(ms)}$, $m=1, 2, \cdots, r$. Since $f_m \in {\bf L}^2\left([0, 1]\right)$, they can be extended as wavelet series  
\be
f_m(t_m)=\sum^{2^{m_{1o}}-1}_{k_1=0}\theta^{(m)}_{j_{mo}, k_m}\psi^{(ms)}_{j_{mo}, k_m}(t_m)+\sum^{\infty}_{j_m=j_{1o}}\sum^{2^{j_m}-1}_{k_m=0}\beta^{(m)}_{j_m,k_m}\psi^{(m)}_{j_m, k_m}(t_m),
\ee
where $\theta^{(m)}_{j_m, k_m}=\int^1_0f_m(t_m)\psi^{(ms)}_{j_m, k_m}(t_m)dt_m$ and  $\beta^{(m)}_{j_m, k_m}=\int^1_0f_m(t_m)\psi^{(m)}_{j_m, k_m}(t_m)dt_m$, $m=1, 2, \cdots, r$. Here, set the lowest resolution levels equal to zero, that is, $j_{mo}=0$, $m=1, 2, \cdots, r$. \\
Thus, for $m=1, 2, \cdots, r$, define the unbiased estimators for the wavelet and scaling coefficients
\be \label{betestm}
\widehat{\beta^{(m)}_{j_m, k_m}}=\frac{1}{N}\sum^N_{i=1}\frac{y({\bf t}_i)}{\Pi^r_{q=1}h_q(t_{qi})}\psi^{(m)}_{j_mk_m}(t_{mi}),
\ee
\be  \label{thetestm}
\widehat{\theta^{(m)}_{j_m, k_m}}=\frac{1}{N}\sum^N_{i=1}\frac{y({\bf t}_i)}{\Pi^r_{q=1}h_q(t_{qi})}\left(\psi^{(ms)}_{j_mk_m}(t_{mi})-2^{-j_m/2}\right).
\ee
Therefore, for $m=1, 2, \cdots, r$, consider the hard thresholding estimators for $f_m$,
\be\label{festm}
\widehat{f}_{mN}(t_m)=\widehat{\theta^{(m)}_{0,0}}\psi^{(ms)}_{0,0}(t_m)+\sum^{J_m-1}_{j_1=j_{1o}}\sum^{2^{j_m}-1}_{k_m=0}\widehat{\beta^{(m)}_{j_m,k_m}}\II\left(|\widehat{\beta^{(m)}_{j_mk_m}}|> \lambda^{(m)}(j_m)\right)\psi^{(m)}_{j_m, k_m}(t_m),
\ee
where, the quantities $J_m$, $j_{mo}$, and $\lambda_m(j_m)$, $m=1,2, \cdots$, are to be determined. In addition, an unbiased estimator for $\beta_o$ is given by
\be \label{consetr}
\widehat{\beta}_o=\frac{1}{N} \sum^N_{i=1} \frac{y(t_i, x_i)}{h_1(t_i)h_2(x_i)}.
\ee
Then, the estimator for the bivariate function $U({\bf t})$ is the sum of its components' own estimators \fr{consetr} and \fr{festm}, as follows
\be \label{uestm}
\widehat{U}_N(t, x)=\widehat{\beta}_o+\sum^r_{m=1}\widehat{f}_{mN}(t_m).
\ee
{\bf Assumption A.6.}\label{A6} $f_m \in {\bf{L^{2}}}\left[0, { 1} \right]$ are bounded above, that is, there exists positive constant $M_m<\infty$, $m=1, 2, \cdots, r$, such that  $f_m(t_m) \leq M_m$ for all  $ t_m \in \left[0, { 1} \right]$.\\
{\bf Assumption A.7.} \label{A7} The probability density functions $h_q$, $q=1, 2, \cdots, r$, are uniformly bounded, that is, on  $\left[0, {1} \right]$ there exist positive constants $m_{q1}$, $m_{q2}$, with $0< m_{q1}\leq m_{q2} < \infty$, $q=1, 2, \cdots, r$, such that $m_{q1} \leq h_q(t) \leq m_{q2}$. \\
{\bf Assumption A.8.}  \label{A8}The functions $f_m$, $m=1, 2, \cdots, r$, are such that $\int^1_0f_m(t_m)dt_m=0$. \\
  {\bf Assumption A.9.} \label{A9}Denote for $m=1, 2, \cdots, r$,
  \beqn  \label{eq10}
 s^{*}_m&=&s_m+1/2 - 1/p,
 \eeqn
  \beqn  \label{s-bar}
 \overline{s}&=& \min\{s_1, s_2, \cdots, s_r\}.
 \eeqn
  The functions $f_m(t)$, $m=1, 2, \cdots, r$, belong to a Besov space. In particular, if $s_{mo} \geq s_m$, their wavelet coefficients $ \beta^{(m)}_{j_mk_m}$ satisfy
\be  \label{eq11}
 B^{s_1, s_2, \cdots, s_r}_{p, q}(A)=\left \{ U \in {\bf L}^2([0, 1]): \left( \sum_{j_m} 2^{(j_ms_m^{*})q}\left (\sum_{k_m}| \beta^{(m)}_{j_m, k_m}|^{p}\right)^{q/{p}}\right )^{1/q} \leq A_m, m=1, 2, \cdots, r\right \}.
\ee
\noindent
{\bf Assumption A.10.} \label{A10} $\sigma({\bf t})$ is bounded below on $[0, 1]^r$, that is, there exist constant $\sigma_{ro}$ such that and $\sigma({\bf t}) >\sigma_{ro}$.
\begin{theorem}\label{th:lowerbds-r} Let  Assumptions {\bf{A.6-A.10}} hold. Then, as $N\rightarrow \infty$, one has
 \be \label{lowerbds-r}
 R({\bf{B^{s_1, s_2, \cdots, s_r}}}(A))\geq C 
{\bar{A}}^2 \left[\frac{1}{{\bar{A}}^2N} \right]^{\frac{2\overline{s}}{2\overline{s}+1}}\vee \frac{1}{N^{\alpha}}, 
 \ee
 where $\bar{A}$ is the radius that corresponds to function with the smallest $s_m$ from amongst the functions $f_1$, $f_2$, $\cdots$, $f_r$. 
 \end{theorem}
\begin{lemma} \label{lem:Var-r}
Let conditions ${\bf{A.6}}$-${\bf{A.8}}$ hold and let $\widehat{\beta^{(m)}_{j_mk_m}}$ and $\widehat{\theta^{(m)}_{j_mk_m}}$ be defined in \fr{betestm}-\fr{thetestm}, for $m=1, 2, \cdots, r$. Then, provided that $\int_{[0, 1]^r}\psi^{(m)}_{j_mk_m}(t_m)\sigma({\bf t})d{\bf t}\neq 0$, as $N \rightarrow \infty$, one has 
\be \label{var-bet-r}
\EE\|\widehat{\beta^{(m)}_{j_mk_m}}-\beta^{(m)}_{j_mk_m}\|^2 \asymp \frac{1}{N}+\frac{1}{N^{\alpha}}.
\ee
\be \label{var-thet-r}
\EE\|\widehat{\theta^{(m)}_{j_mk_m}}-\theta^{(m)}_{j_mk_m}\|^2 \asymp \frac{1}{N}+ \frac{1}{N^{\alpha}}.
\ee
\be \label{var2-bet-r}
\EE\|\widehat{\beta^{(m)}_{j_mk_m}}-\beta^{(m)}_{j_mk_m}\|^4 \asymp \frac{1}{N^2}.
\ee
In addition, for the estimator $\widehat{\beta}_o$ in \fr{consetr}, under condition \fr{longm} as $N \rightarrow \infty$, one has 
\be \label{var-beto}
\EE|\widehat{\beta}_{o}-\beta_{o}|^2 \asymp \frac{\sigma^2}{N^{\alpha}}.
\ee
\end{lemma}
Choose for estimator \fr{festm} the thresholds $\lambda^{(m)}(j_m)$, $m=1, 2, \cdots, r$, such that 
\be \label{thresh21-r}
\lambda^{(m)}(j_m)=\left\{ \begin{array}{ll}  
\gamma_m\frac{{\ln(N)}}{\sqrt{N}}, \ & \mbox{if}\ \int_{[0, 1]^r}\psi^{(m)}_{j_mk_m}(t_m)\sigma({\bf t})d{\bf t}=0, \\
\gamma_{m1}\sqrt{\frac{\ln{(N)}}{N^{\alpha}}}\ & \mbox{if}\ otherwise.
  \end{array} \right.
  \ee
 In addition, the highest resolution levels $J_m$, $m=1, 2, \cdots, r$, should be chosen such that
\be \label{Jr}
2^{J_m}\sim \frac{N}{\ln(N)}.
\ee 
\begin{lemma} \label{lem:Lardev-r}
Let $\widehat{\beta^{(m)}_{j_mk_m}}$, $m=1, 2, \cdots, r$, be defined in \fr{betestm}, and let assumptions ${\bf{A.6}}$-${\bf{A.8}}$ and condition \fr{longm} hold. Then, as $N \rightarrow \infty$, one has 
\be  \label{Largdev-r}
\Pr \left(| \widehat{\beta^{(m)}_{j_mk_m}}-\beta^{(m)}_{j_mk_m} |> 1/2\lambda^{(m)}({j_m})\right)\asymp \left[\frac{1}{N}\right]^{2\tau},
\ee
where $\tau$ is a positive parameter that is large enough. 
\end{lemma}
\begin{theorem} \label{th:upperbds-r}
Let  $\widehat{U}_{N}(t, x)$ be the estimator defined in \fr{uestm} with its components $\widehat{\beta}_o$ and $\widehat{f}_{mN}(t_m)$, $m=1, 2, \cdots, r$, given in \fr{consetr} and \fr{festm}, respectively. Suppose assumptions ${\bf{A.6}}$-${\bf{A.10}}$ hold. Then, under condition \fr{longm}, if $\tau$ is large enough, as $N \rightarrow \infty$, one has
 \be \label{upperbds-lm}
 \sup_{U\in B^{s_1, s_2, \cdots, s_r}(A)} \EE\|\widehat{U}_N-U\|^2\leq C 
{\bar{A}}^2 \left[\frac{\ln(N)}{{\bar{A}}^2N} \right]^{\frac{2\overline{s}}{2\overline{s}+1}}\vee \frac{1}{N^{\alpha}}.
 \ee
\end{theorem}
\begin{remark}	  
As you can see from {\bf Theorem \ref{th:lowerbds-r}} and {\bf Theorem \ref{th:upperbds-r}}, the convergence rates depend only on the lowest smoothness parameter $\overline{s}= \min\{s_1, s_2, \cdots, s_r\}$ from amongst the univariate functions $f_1$, $f_2$, $\cdots$, $f_r$, or the long-memory parameter $\alpha$. This implies that our estimator does not suffer from the curse of dimensionality. This behavior is quite common in the literature. For instance, this behavior is similar to Grez and Vidakovic~(2018) when their joint design density is uniform over the torus $[0, 1]^r$ and our $\alpha > \frac{2\overline{s}}{2\overline{s}+1}$. 
  \end{remark}
\begin{remark}	
Our design densities $h_m$, $m=1, 2, \cdots, r$, are assumed to me known, but in practice this may not be the case. These functions can be estimated from the data and their empirical counterparts may be used in formulas \fr{betestm}, \fr{thetestm} and \fr{consetr}, and the interested reader may refer to Pensky and Vidakovic~(2001). Providing completely data driven estimates for $h_m$, $m=1, 2, \cdots, r$, is beyond the scope of this work and we assume that these functions are known. 
 \end{remark}

 \section{Proofs }
In order to prove Theorem \ref{th:lowerbds}, we use the following lemma  
\begin{lemma} (Lemma $A.1$ of Bunea et al.~(2007)) 
\label{lem:Bunea} 
Let $\Te$ be a set of functions of cardinality $\card(\Te)\geq 2$ such that\\
(i) $\|f-g\|_p^p \geq 4\delta^p, \ for\  f, g \in \Te, \ f \neq g, $\\
(ii) the Kullback divergences $K(P_f, P_g)$ between the measures $P_f$ and $P_g$ 
satisfy the inequality $K(P_f, P_g) \leq \log(\card(\Te))/16,\ for\ f,\ g \in \Te$.\\
Then, for some absolute positive constant $C_1$, one has 
$$
 \inf_{f_n}\sup_{f\in \Te} \EE_f\|f_n-f\|_p^p \geq C_1 \delta^p,
$$
where $\inf_{f_n}$ denotes the infimum over all estimators.
\end{lemma}
 {\bf Proof of Theorem \ref{th:lowerbds} and Theorem \ref{th:lowerbds-lm}}. According to Yang~(2001), the minimax risk of estimating a regression function under random design and long-memory errors is the sum of two terms; the minimax rate of the same class but under $i.i.d.$ errors, and the rate of estimating the mean of the regression function under correlated errors. So, we will use this result in our derivation of the lower bounds. Also, since our regression function is additive, we will consider two cases; $U(t, x)\equiv f(t)$ and $U(t, x)\equiv g(x)$. \\
  Case 1: $U(t, x)\equiv \beta_o+f(t)$. Let $\omega$ be the vector with components $\omega_k\in \{ 0, 1\}$, $k=0, 1, \cdots, 2^j-1$, and denote the set of all possible values of $\omega$ by $\Omega$. Let $f_{\omega}$ be the functions of the form
\be  \label{f_j}
f_{\omega}(t)=\rho_j\sum^{2^J-1}_{k=1}\omega_k\psi_{jk}(t),\ \ \omega_l \in \{ 0, 1\}.
\ee
Observe that $\omega$ has $\Gamma=2^j$ components and therefore $\Omega$ will have cardinality $\card(\Omega)=2^{\Gamma}$. By \fr{eq11}, it is easy to verify that  $f_{\omega}(t) \in B^{s_1, s_2}(A)$ with the choice $\rho^2_{j}=A_1^22^{-j(2s_1+1)}$. Take ${f}_{\tilde{\omega}}$ of the form of \fr{f_j} but with $\tilde{\omega}_l\in \{ 0, 1\}$, then applying Varshamov-Gilbert Lemma (\cite{tsybakov}, p 104), the $L^2$-norm of the difference is
\be
\| f_{\omega}(t)-\tilde{f}_{{\omega}}(t)\|^2 \geq  \frac{2^{j_1}\rho^2_{j_1}}{8}.
\ee
To prove {\bf Theorem \ref{th:lowerbds}}, define the quantities $\pmb{h}_i=\beta_o+ f_{\omega}(t_i)+ \sigma\mu_i$, $i=1, 2, \cdots, N$, where $\{\mu_i\}_{i\geq 1}$ is stationary Gaussian sequence that satisfies structure \fr{longm}. Let $P_{f_{\omega}}$ be the probability law of the process $\pmb{h}_i$ under the hypothesis $f_{\omega}$ defined in \fr{f_j}. Then, by Assumptions ${\bf A.2}$ and ${\bf A.3}$, the Kullback divergence can be written as 
\beqn
K(P_{f_{\omega}}, P_{\tilde{f}_{\omega}})&\leq& \frac{1}{2\sigma^2}\left\{\EE\left[f_{\omega}(t_i)-f_{\tilde{\omega}}(t_i)\right]^2\sum^N_{i=1}\xi^{-1}_{ii}+\left[\EE\left(f_{\omega}(t_i)-f_{\tilde{\omega}}(t_i)\right)\right]^2\sum^N_{i\neq i'}\xi^{-1}_{ii'} \right\}\nonumber\\
& \leq & \frac{1}{2\sigma^2}\rho^2_{j}\left\{\sum^{2^j-1}_{k=1}|\omega_k-\tilde{\omega}_k|^2\int^1_0\psi^2_{jk}(t)h_1(t)dt\sum^N_{i=1}\xi^{-1}_{ii}+\sum^{2^j-1}_{k=1}|\omega_k-\tilde{\omega}_k|^2\left[\int^1_0\psi_{jk}(t)h_1(t)dt\right]^2\sum^N_{i\neq i'}\xi^{-1}_{ii'} \right\}\nonumber\\
&\leq& \frac{m_{12}}{2\sigma^2}\rho^2_j2^j\left\{ Trace\left(\Sigma^{-1}_N\right)+{\bf1}^T\Sigma^{-1}_N{\bf1}\right\}\nonumber\\
&\leq & \frac{m_{12}}{2\sigma^2}\rho^2_j2^j\left\{ N+N^{\alpha} \right\}\leq \frac{m_{12}N}{\sigma^2}\rho^2_j2^j,
\eeqn
where $\Sigma_N$ is the $N\times N$ covariance matrix of $\{ \varepsilon_i\}_{i\geq1}$, $\xi^{-1}_{ii'}$ are the elements of the inverse of $\Sigma_N$, and ${\bf 1}$ is the N-dimensional vector whose elements that are all equal to 1. Now, to apply {\bf Lemma \ref{lem:Bunea}}, choose 
\be \label{klblam}
A_1^2 2^{-j(2s+1)}2^j\frac{m_{12}N}{\sigma^2}\leq \pi_02^j,
\ee
that is,
\be \label{Lop}
2^j = C\left[\frac{A_1^2N}{\sigma^2}\right]^{\frac{1}{2s_1 +1}}.
\ee
Hence, the lower bounds are
\be \label{cas1}
\delta= CA_1^2\left[\frac{\sigma^2}{A_1^2N}\right]^{\frac{2s_1}{2s_1 +1}}\vee \frac{1}{N^{\alpha}}.
\ee
 2: $U(t, x)\equiv \beta_o+g(x)$. We follow the same procedure as in Case 1 to obtain the lower bounds
\be \label{cas2}
\delta= CA_2^2\left[\frac{\sigma^2}{A_2^2N}\right]^{\frac{2s_2}{2s_2 +1}}\vee \frac{1}{N^{\alpha}}. 
\ee
To complete the proof, we choose the highest between \fr{cas1} and \fr{cas2}.  $\Box$\\
To prove {\bf Theorem \ref{th:lowerbds-lm}}, we consider two cases in the spirit of the proof of {\bf Theorem \ref{th:lowerbds}}.   Case 1: $U(t, x)\equiv \beta_o+f(t)$. Based on the same test functions as the previous proof, define the $\pmb{h}_i=\beta_o+ f_{\omega}(t_i)+ \sigma(t_i, x_i)\mu_i$, $i=1, 2, \cdots, N$, where $\{\mu_i\}_{i\geq 1}$ is stationary Gaussian sequence that satisfies structure \fr{longm}. Let $P_{f_{\omega}}$ be the probability law of the process $\pmb{h}_i$ under the hypothesis $f_{\omega}$ defined in \fr{f_j}. Then, by Assumptions ${\bf A.2}$ and ${\bf A.3}$, the Kullback divergence can be written as 
\beqn
K(P_{f_{\omega}}, P_{\tilde{f}_{\omega}})&\leq& \frac{1}{2}\left\{\EE\left[\frac{f_{\omega}(t_i)-f_{\tilde{\omega}}(t_i)}{\sigma^2(t_i, x_i)}\right]^2\sum^N_{i=1}\xi^{-1}_{ii}+\left[\EE\left(\frac{f_{\omega}(t_i)-f_{\tilde{\omega}}(t_i)}{\sigma(t_i, x_i)}\right)\right]^2\sum^N_{i\neq i'}\xi^{-1}_{ii'} \right\}\nonumber\\
& \leq & \frac{1}{2}\rho^2_{j}\left\{\sum^{2^j-1}_{k=1}|\omega_k-\tilde{\omega}_k|^2\int^1_0\int^1_0\frac{\psi^2_{jk}(t)h_1(t)}{\sigma^{2}(t, x)}dth_2(x)dx\sum^N_{i=1}\xi^{-1}_{ii}\right\}\nonumber\\
&+&\left\{\sum^{2^j-1}_{k=1}|\omega_k-\tilde{\omega}_k|^2\left[\int^1_0\int^1_0\frac{\psi_{jk}(t)h_1(t)}{\sigma(t, x)}dth_2(x)dx\right]^2\sum^N_{i\neq i'}\xi^{-1}_{ii'} \right\}\nonumber\\
&\leq& \frac{m_{12}m_{22}}{2}\|\sigma^{-2}_{tx}\|_{\infty}\rho^2_j2^j\left\{ Trace\left(\Sigma^{-1}_N\right)+{\bf1}^T\Sigma^{-1}_N{\bf1}\right\}\nonumber\\
&\leq & \frac{m_{12}m_{22}}{2}\|\sigma^{-2}_{tx}\|_{\infty}\rho^2_j2^j\left\{ N+N^{\alpha} \right\}\leq {m_{12}m_{22}N}\|\sigma^{-2}_{tx}\|_{\infty}\rho^2_j2^j. 
\eeqn
Now, to apply {\bf Lemma \ref{lem:Bunea}}, choose 
\be \label{klblam}
A_1^2 2^{-j(2s_1+1)}2^j{m_{12}m_{22}N}{\|\sigma^{-2}_{tx}\|_{\infty}}\leq \pi_02^j,
\ee
that is,
\be \label{Lop}
2^j = C\left[\frac{A_1^2N}{\sigma^2}\right]^{\frac{1}{2s_1 +1}}.
\ee
Hence, the lower bounds are
\be \label{cas1}
\delta= CA_1^2\left[\frac{\sigma^2}{A_1^2N}\right]^{\frac{2s_1}{2s_1 +1}}\vee \frac{1}{N^{\alpha}}.
\ee
 2: $U(t, x)\equiv \beta_o+g(x)$. We follow the same procedure as in Case 1 to obtain the lower bounds
\be \label{cas2}
\delta= CA_2^2\left[\frac{\sigma^2}{A_2^2N}\right]^{\frac{2s_2}{2s_2 +1}}\vee \frac{1}{N^{\alpha}}. 
\ee
To complete the proof, we choose the highest between \fr{cas1} and \fr{cas2}.  $\Box$\\
{\bf Proof of Lemma \ref{lem:Var} and \ref{lem:Var-lm}.} Notice that with $\widehat{\beta^{(1)}_{j_1k_1}}$  defined in \fr{betest1}, one has 
\beqn \label{alef}
\widehat{\beta^{(1)}_{j_ik_i}}-\beta^{(1)}_{j_ik_i}&=& \frac{1}{N}\sum^N_{i=1}\left[\frac{f(t_i)\psi_{j_1k_1}(t_i)}{h_1(t_i)h_2(x_i)}-\EE\left[\frac{f(t_i)\psi_{j_1k_1}(t_i)}{h_1(t_i)h_2(x_i)}\right]\right]\nonumber\\
 &+& \frac{1}{N}\sum^N_{i=1}\left[\frac{g(x_i)\psi_{j_1k_1}(t_i)}{h_1(t_i)h_2(x_i)}-\EE\left[\frac{g(x_i)\psi_{j_1k_1}(t_i)}{h_1(t_i)h_2(x_i)}\right]\right]\nonumber\\
 &+& \frac{1}{N}\sum^N_{i=1}\frac{\psi_{j_1k_1}(t_i)\sigma(t_i, x_i)}{h_1(t_i)h_2(x_i)}\varepsilon_i=\left(\Delta_{11}+\Delta_{12}\right)+\Delta_2=i.i.d\ term+ martingale/ long\ memory\ term.\nonumber\\
 \eeqn
Let $\mathcal{F}_i=\sigma\left( \xi_i, (t_i, x_i), \xi_{i-1}, (t_{i-1}, x_{i-1}), \cdots \right)$ and let $\varepsilon_{i, i-1}=\varepsilon_i-\xi_i$. Note that $\varepsilon_{i, i-1}$ is $\mathcal{F}_{i-1}$-measurable and $(\xi_i, \{t_i, x_i\})$ is independent of $\mathcal{F}_{i-1}$. Therefore, 
\be
\EE\left[ \frac{\psi_{j_1k_1}(t_i)\sigma(t_i, x_i)}{h_1(t_i)h_2(x_i)}\varepsilon_i\mid \mathcal{F}_{i-1}\right]=\varepsilon_{i, i-1}\EE\left[\frac{\psi_{j_1k_1}(t_i)\sigma(t_i, x_i)}{h_1(t_i)h_2(x_i)}\right].
\ee
Thus, fort case $\sigma(t, x)\equiv \sigma$, , a positive constant, $\Delta_2$ can be partitioned as follows
\beqn
\Delta_2&=&\frac{1}{N}\sum^N_{i=1}\left[\frac{\psi_{j_1k_1}(t_i)\sigma(t_i, x_i)}{h_1(t_i)h_2(x_i)}\varepsilon_i-\EE\left[ \frac{\psi_{j_1k_1}(t_i)\sigma(t_i, x_i)}{h_1(t_i)h_2(x_i)}\varepsilon_i\mid \mathcal{F}_{i-1}\right]\right]\nonumber\\
&+& \frac{1}{N}\EE\left[\frac{\psi_{j_1k_1}(t_1)\sigma(t_i, x_i)}{h_1(t_1)h_2(x_1)}\right]\sum^N_{i=1}\varepsilon_{i, i-1}\nonumber\\
&=& \frac{\sigma}{N}\sum^N_{i=1}\left[\frac{\psi_{j_1k_1}(t_i)}{h_1(t_i)h_2(x_i)}\varepsilon_i-\EE\left[ \frac{\psi_{j_1k_1}(t_i)}{h_1(t_i)h_2(x_i)}\varepsilon_i\mid \mathcal{F}_{i-1}\right]\right],
\label{delta2}
\eeqn
since wavelets are orthogonal with respect to constants, 
\beqns
\EE\left[\frac{\psi_{j_1k_1}(t_1)}{h_1(t_1)h_2(x_1)}\right]=\int^1_0\psi_{j_1k_1}(t)dt\int^1_0dx=0. 
\eeqns
Notice that $\Delta_{11}$ and $\Delta_{12}$ are the sums of zero-mean independent and identically distributed components. Therefore, taking expectation of the square of $\Delta_{11}$ yields
\beqn
\EE\left[\Delta^2_{11}\right]&=&\frac{1}{N^2}\sum^N_{i=1}\EE\left[\frac{f(t_i)\psi_{j_1k_1}(t_i)}{h_1(t_i)h_2(x_i)}-\EE\left[\frac{f(t_i)\psi_{j_1k_1}(t_i)}{h_1(t_i)h_2(x_i)}\right]\right]^2\nonumber\\
&\leq & \frac{1}{N} \max_{t}|f(t)|^2\EE \left[\frac{\psi_{j_1k_1}(t_i)}{h_1(t_i)h_2(x_i)}-\EE\left[\frac{\psi_{j_1k_1}(t_i)}{h_1(t_i)h_2(x_i)}\right]\right]\nonumber\\
&=& \frac{1}{N} \max_{t}|f(t)|^2\int^1_0\int^1_0\psi^2_{j_1k_1}(t)h_1^{-1}(t)h_2^{-1}(x)dtdx\nonumber\\
&\leq &\frac{1}{N} \max_{t}|f(t)|^2\frac{1}{m_{11}}\frac{1}{m_{21}}. \label{delta11}
\eeqn
The same way, one can show that 
\be \label{delta12}
\EE\left[\Delta^2_{12}\right]\leq \frac{1}{N} \max_{x}|g(x)|^2\frac{1}{m_{11}}\frac{1}{m_{21}}.
\ee
Also, it can be easily shown that for $2^{j_1}<N$, one has
\beqn
\EE\left[\Delta^4_{11}\right]\leq\frac{2}{N^2} \max_{t}|f(t)|^4\frac{1}{m^3_{11}}\frac{1}{m^3_{21}}.\label{delta211}\\
\EE\left[\Delta^4_{12}\right]\leq\frac{2}{N^2} \max_{x}|g(x)|^4\frac{1}{m^3_{11}}\frac{1}{m^3_{21}}.\label{delta212}
\eeqn
To evaluate \fr{delta2}, notice that 
\beqns
\frac{N}{\sigma}\Delta_2=\sum^N_{i=1}d_i 
\eeqns
is a martingale with 
\beqn
d_i&=&\frac{\psi_{j_1k_1}(t_i)}{h_1(t_i)h_2(x_i)}\varepsilon_i-\EE\left[\frac{\psi_{j_1k_1}(t_i)}{h_1(t_i)h_2(x_i)}\varepsilon_i\mid \mathcal{F}_{i-1}\right]\nonumber\\
&=& \left[\frac{\psi_{j_1k_1}(t_i)}{h_1(t_i)h_2(x_i)}-\EE\left[\frac{\psi_{j_1k_1}(t_i)}{h_1(t_i)h_2(x_i)}\right]\right]\varepsilon_{i, i-1}+\xi_i \frac{\psi_{j_1k_1}(t_i)}{h_1(t_i)h_2(x_i)}\nonumber\\
&=&\frac{\psi_{j_1k_1}(t_i)}{h_1(t_i)h_2(x_i)}\varepsilon_{i, i-1}+\xi_i \frac{\psi_{j_1k_1}(t_i)}{h_1(t_i)h_2(x_i)}. \label{di}
\eeqn
Therefore, taking the $p^{th}$, for $p\geq 2$, central moment of \fr{di}
\beqn
\EE|d_i|^p&\leq& 2^{p-1}\left(\EE|\varepsilon_{i, i-1}|^p\EE\left|\frac{\psi_{j_1k_1}(t_i)}{h_1(t_i)h_2(x_i)}\right|^p+\EE|\xi_i|^p\EE\left|\frac{\psi_{j_1k_1}(t_i)}{h_1(t_i)h_2(x_i)}\right|^p\right)\nonumber\\
&\leq & 2^{p-1}\left(\EE|\varepsilon_{i, i-1}|^p\frac{2^{j_1(p/2-1)}}{m^{p-1}_{11}m^{p-1}_{21}}+\EE|\xi_i|^p\frac{2^{j_1(p/2-1)}}{m^{p-1}_{11}m^{p-1}_{21}} \right). \label{edp}
\eeqn
In addition,
\beqn
\EE\left[d_i^2\mid \mathcal{F}_{i-1}\right]=\EE\left[\frac{\psi_{j_1k_1}(t)}{h_1(t)h_2(x)}\right]^2\varepsilon_{i, i-1}^2+\EE\left[\frac{\psi_{j_1k_1}(t)}{h_1(t)h_2(x)}\right]^2\EE\left[\xi_1\right]^2. \label{edcond}
\eeqn
Therefore, By Rosenthal's inequality for martingales and results \fr{edp}(with $p=2$) and \fr{edcond}, one has
\beqn \label{delt2}
\EE\left[ \Delta_2\right]^2\leq \frac{\sigma^2}{N^2}\EE\left[\sum^N_{i=1}\EE\left[d_i^2\mid \mathcal{F}_{i-1}\right]\right]+\frac{\sigma^2}{N^2}\sum^N_{i=1}\EE|d_i|^2\leq \frac{3\sigma^2}{N}\frac{1}{m_{11}m_{21}}\left[\EE|\varepsilon_{1, 0}|^2+\EE\left[\xi_1\right]^2\right]. 
\eeqn
Hence, combining \fr{delta11}, \fr{delta12} and \fr{delt2} completes the proof of \fr{var-bias}. The proof of $\EE\|\widehat{\beta^{(2)}_{j_2k_2}}-\beta^{(2)}_{j_2k_2}\|^2$ and \fr{var-thet} 
 can be derived the same way. $\Box$\\
 To show \fr{var2-bet} notice that by \fr{edp} and \fr{edcond}, and for $2^{j_1}<N$, one has
   \beqn \label{delt24}
 \EE\left[ \Delta_2\right]^4&\leq& \frac{\sigma^4}{N^4}\EE\left[\sum^N_{i=1}\EE\left[d_i^2\mid \mathcal{F}_{i-1}\right]\right]^2+\frac{\sigma^4}{N^4}\sum^N_{i=1}\EE|d_i|^4\nonumber\\
 &\leq& \frac{\sigma^4}{N^4}\left\{\EE\left[\frac{1}{m_{11}m_{21}}\int^1_0\psi_{j_1k_1}^2(t)dt\sum^N_{i=1}\left[\varepsilon_{i, i-1}^2+\EE(\xi_1^2)\right]\right]^2+\frac{8N}{m_{11}^3m_{21}^3}\int^1_0\psi^4_{j_1k_1}(t)dt\left[\EE|\varepsilon_{1, 0}|^4+\EE\left[\xi_1\right]^4\right]\right\}\nonumber\\
 & \leq & \frac{C\sigma^4}{N^2}. 
   \eeqn
 To complete the proof of  \fr{var2-bet} with $i=1$, use \fr{delta11}-\fr{delta212}, \fr{delt2} and \fr{delt24}, notice that 
 \beqn
 \EE\|\widehat{\beta^{(1)}_{j_1k_1}}-\beta^{(1)}_{j_1k_1}\|^4&=&\EE\left[\Delta_{11}\right]^4+\EE\left[\Delta_{12}\right]^4+\EE\left[\Delta_{2}\right]^4+6\EE\left[\Delta_{11}\right]^2\EE\left[\Delta_{12}\right]^2+6\EE\left[\Delta_{11}\right]^2\EE\left[\Delta_{2}\right]^2\nonumber\\
 &+& 6\EE\left[\Delta_{12}\right]^2\EE\left[\Delta_{2}\right]^2
      \eeqn 
The same way we can prove \fr{var2-bet} with $i=2$. \\
For case $\sigma(t, x)$, a positive non-constant bivariate function, $\Delta_2$ can be partitioned as follows
\beqn
\Delta_2&=&\frac{1}{N}\sum^N_{i=1}\left[\frac{\psi_{j_1k_1}(t_i)\sigma(t_i, x_i)}{h_1(t_i)h_2(x_i)}\varepsilon_i-\EE\left[ \frac{\psi_{j_1k_1}(t_i)\sigma(t_i, x_i)}{h_1(t_i)h_2(x_i)}\varepsilon_i\mid \mathcal{F}_{i-1}\right]\right]\nonumber\\
&+& \frac{1}{N}\EE\left[\frac{\psi_{j_1k_1}(t_1)\sigma(t_1, x_1)}{h_1(t_1)h_2(x_1)}\right]\sum^N_{i=1}\varepsilon_{i, i-1}\nonumber\\
&=&\Delta_{21}+\Delta_{22}=martingale + long\ memory,
\label{delta2}
\eeqn
where
\beqns
\EE\left[\frac{\psi_{j_1k_1}(t_1)\sigma(t_1, x_1)}{h_1(t_1)h_2(x_1)}\right]=\int^1_0\int^1_0\psi_{j_1k_1}(t)\sigma(t, x)dtdx\neq0. 
\eeqns
The i.i.d. terms will be the same as before. For the long-memory term, by condition \fr{longm}, we have 
\beqn \label{delt'22}
\EE\left[\Delta_{22}\right]^2&=&\frac{1}{N^2}\left(\EE\left[\frac{\psi_{j_1k_1}(t_1)\sigma(t_1, x_1)}{h_1(t_1)h_2(x_1)}\right]\right)^2\Var\left[\sum^N_{i=1}\varepsilon_{i, i-1}\right]\leq \EE\left[\frac{\psi_{j_1k_1}(t_1)\sigma(t_1, x_1)}{h_1(t_1)h_2(x_1)}\right]^2 c_{\alpha}N^{-\alpha}\nonumber\\
&\leq & N^{-\alpha}\frac{c_{\alpha}}{m_{11}m_{21}}\int^1_0\int^1_0\psi_{j_1k_1}^2(t)\sigma^2(t, x)dtdx\leq N^{-\alpha}\frac{c_{\alpha}}{m_{11}m_{21}}\|\sigma\|_{\infty}^2.
\eeqn
As for the martingale term, by Rosenthal's inequality for martingales, it can be shown that 
\beqn \label{delt'2}
\EE\left[ \Delta_{21}\right]^2&\leq&  \frac{3}{N}\frac{1}{m_{11}m_{21}}\left[\EE|\varepsilon_{1, 0}|^2+\EE\left[\xi_1\right]^2\right]\int^1_0\int^1_0\psi_{j_1k_1}^2(t)\sigma^2(t, x)dtdx\nonumber\\
&\leq&  \frac{3}{N}\frac{1}{m_{11}m_{21}}\left[\EE|\varepsilon_{1, 0}|^2+\EE\left[\xi_1\right]^2\right]\|\sigma\|_{\infty}^2. 
\eeqn
Consequently, combining \fr{delt'22}, \fr{delt'2}, \fr{delta11} and \fr{delta12} completes the proof for $i=1$. $i=2$ could be dealt with in the same fashion. \\
Now, to prove $(20)$ and $(30)$, notice that
\beqn
\widehat{\beta}_o-\beta_o&=&\frac{\beta_o}{N}\sum^N_{i=1}\left[\frac{1}{h_1(t_i)h_2(x_i)}-1\right]+\frac{1}{N}\sum^N_{i=1}\frac{f(t_i)+g(x_i)}{h_1(t_i)h_2(x_i)}+\frac{1}{N}\sum^N_{i=1}\frac{\sigma(t_i, x_i)}{h_1(t_i)h_2(x_i)}\varepsilon_{i}\nonumber\\
&=&\frac{\beta_o}{N}\sum^N_{i=1}\left[\frac{1}{h_1(t_i)h_2(x_i)}-1\right]+\frac{1}{N}\sum^N_{i=1}\frac{f(t_i)+g(x_i)}{h_1(t_i)h_2(x_i)}\nonumber\\
&+&\frac{1}{N}\sum^N_{i=1}\left[\frac{\sigma(t_i, x_i)}{h_1(t_i)h_2(x_i)}\varepsilon_i-\EE\left(\frac{\sigma(t_i, x_i)}{h_1(t_i)h_2(x_i)}\varepsilon_i\mid \mathcal{F}_{i-1}\right)\right] + \frac{1}{N}\EE\left[\frac{\sigma(t_1, x_1)}{h_1(t_1)h_2(x_2)}\right]\sum^N_{i=1}\varepsilon_{i, i-1}\nonumber\\
&=&\left(\Delta'_1+\Delta'_2\right)+\Delta'_{31}+\Delta'_{32}= i.i.d\ term+martingale\ term\ +long-memory\ term. \nonumber\\
\eeqn
The i.i.d. term and the martingale term will be dealt with the same way as the proof of \fr{var-bet}, so we will only evaluate the long-memory term. Notice that under both cases of $\sigma(t, x)$, one has
\beqn
\sigma_s=\EE\left[\frac{\sigma(t, x)}{h_1(t)h_2(x)}\right]=\int^1_0\int^1_0\sigma(t, x)dtdx\neq 0. 
\eeqn
Therefore, by condition \fr{longm}
\beqn
\EE\left[\Delta'_{32}\right]^2=\frac{1}{N^2}\left[\int^1_0\int^1_0\sigma(t, x)dtdx\right]^2\Var\left(\sum^N_{i=1}\varepsilon_{i, i-1}\right)\sim \pi_{\alpha}\frac{\sigma_s^2}{N^{\alpha}}.
\eeqn
$\Box$\\
  {\bf Proof of Lemma \ref{lem:Lardev} and \ref{lem:Lardev-lm}.}  In order to prove \fr{Largdev}, we make use of Bernstein inequality version for independent and identically distributed random variables and the martingale version. 
    \begin{lemma}(Bernstein Inequality: i.i.d. Version). \label{lem:bernineqiid}
Let $Y_i$, $i=1, 2, \cdots, N$,  be independent and identically distributed random variables with mean zero and finite variance $\sigma^2$, with $\|Y_i\| \leq \|Y\|_{\infty} < \infty$. Then,  
\be \label{prob-biid}
\Pr\left(\left|N^{-1}\sum^N_{i=1}Y_i\right|>z\right) \leq 2 \exp\left\{-\frac{Nz^2}{2(\sigma^2+\|Y\|_{\infty}z/3)}\right\}.
\ee
\end{lemma}
 \begin{lemma}(Bernstein Inequality: Martingale Version). (Dzhaparidze and Van Zanten~(2001)). \label{lem:bernineqmar}
Let $\left(d_i, \mathcal{F}_i\right)$, $i=1, 2, \cdots, N$,  be a martingale difference sequence.  Then,  
\be \label{prob-bmar}
\Pr\left(\left|\sum^N_{i=1}d_i\right|>z, \sum^N_{i=1}d_i^2\II(|d_i|>x)+\sum^N_{i=1}\EE\left[d_i^2\mid \mathcal{F}_{i-1}\right]\leq L\right) \leq 2 \exp\left\{-\frac{z^2}{2(L+xz/3)}\right\}.
\ee
\end{lemma}
Now, recall decomposition \fr{alef}. Therefore, for $i=1$ and $\sigma(t, x)\equiv \sigma$, the probability
\beqn
\Pr \left(| \widehat{\beta^{(i)}_{j_ik_i}}-\beta^{(i)}_{j_ik_i} |> 1/2\lambda({j_i})\right)&\leq& \Pr\left(\left|\Delta_{11}\right|>1/6\lambda^{(1)}(j_1)\right)+\Pr\left(\left|\Delta_{12}\right|>1/6\lambda^{(1)}(j_1)\right)\nonumber\\
&+& \Pr\left(\left|\Delta_{2}\right|>1/6\lambda^{(1)}(j_1)\right)=P_1+P_2+P_3. 
\eeqn
For $P_1$ and $P_2$ we apply result \fr{prob-biid}, since $\Delta^{(1)}_{11}$ and $\Delta^{(1)}_{12}$ are zero mean with variances given by \fr{delta11} and \fr{delta12}, respectively, and
\beqn
\|\Delta^{(1)}_{11}\|_{\infty}\leq \frac{1}{N}\|f\|_{\infty}\|\psi_{j_1k_1}\|_{\infty}\frac{1}{m_{11}m_{21}}, \ \ \ \|\Delta^{(1)}_{12}\|_{\infty}\leq \frac{1}{N}\|g\|_{\infty}\|\psi_{j_1k_1}\|_{\infty}\frac{1}{m_{11}m_{21}}.
\eeqn 
Indeed, by \fr{prob-biid}, one has
\beqn \label{p1}
P_1\leq 2N^{-\frac{\gamma^2m_{11}m_{21}}{144\|f\|_{\infty}}},
\eeqn
and 
\beqn \label{p2}
P_2\leq 2N^{-\frac{\gamma^2m_{11}m_{21}}{144\|g\|_{\infty}}}.
\eeqn
Now, for $P_3$, denote $\aleph(N, x)=\sum^N_{i=1}d_i^2\II(|d_i|>x)+\sum^N_{i=1}\EE\left[d_i^2\mid \mathcal{F}_{i-1}\right]$ and notice that
\beqn \label{p3}
P_3&=&\Pr\left(\frac{1}{N}\left|\sum^N_{i=1}d_i\right|>\frac{1}{6\sigma}\lambda^{(1)}(j_1) \right)\nonumber\\
&\leq& \Pr\left( \frac{1}{N}\left|\sum^N_{i=1}d_i\right|>\frac{1}{6\sigma}\lambda^{(1)}(j_1), \aleph(N, x)\leq L \right)+\Pr\left(\aleph(N, x)>L\right).
\eeqn
Thus, by result \fr{prob-bmar}, the first term of \fr{p3} becomes
\beqn
\Pr\left( \frac{1}{N}\left|\sum^N_{i=1}d_i\right|>\frac{1}{6\sigma}\lambda^{(1)}(j_1), \aleph(N, x)\leq L \right)\leq 2 \exp\left\{-\frac{N^2{1/(36\sigma^2)}\lambda^{(1)}(j_1)^2}{2(L+xN{\lambda^{(1)}(j_1)/18\sigma})}\right\}.
\eeqn
As for the second term of \fr{p3}, it can be decomposed into two terms as follows
\beqn\label{p31}
\Pr\left(\aleph(N, x)>L\right)\leq \Pr\left(\sum^N_{i=1}d_i^2\II(|d_i|>x)> L/2\right)+\Pr\left(\sum^N_{i=1}\EE\left[d_i^2\mid \mathcal{F}_{i-1}\right]>L/2\right)=P_{321}+P_{322}.
\eeqn
For $P_{322}$, notice that by choosing $L=2N\left[D\ln(N)+\frac{1}{m_{11}m_{21}}\EE(\xi^2_1)\right]$ and keeping in mind that $\varepsilon_{i, i-1}$ is Gaussian, one has 
\beqn \label{p322}
P_{322}= \Pr\left(\sum^N_{i=1}\EE\left[d_i^2\mid \mathcal{F}_{i-1}\right]>L/2\right)&\leq& \Pr\left(\frac{1}{m_{11}m_{21}}\left[N\EE(\xi^2_1)+\sum^N_{i=1}\varepsilon^2_{i, i-1}\right]> L/2\right).\nonumber\\
&\leq & N\Pr\left(|\varepsilon_{1, 0}|> \sqrt{Dm_{11}m_{21}\ln(N)}\right)\nonumber\\
&\leq & 2N\exp\left\{-\frac{Dm_{11}m_{21}\ln(N)}{2\EE[\varepsilon_{1,0}]^2}\right\}
\eeqn
As for $P_{321}$, notice that $\|\psi_{j_1k_1}\|_{\infty}\leq 2^{j_1}\|\psi\|_{\infty}\leq |\psi\|_{\infty}\frac{N}{\ln(N)}$. Since $\varepsilon_{i, i-1}$ and $\xi_i$ are Gaussian and independent from each other, taking $x^2=\pi_o\sqrt{N}$, yields
\beqn \label{p321}
P_{321}&=& \Pr\left(\sum^N_{i=1}d_i^2\II(|d_i|>x)> L/2\right)\leq   \Pr\left(\sum^N_{i=1}d_i^2\II(|d_i|>x)> DN\ln(N)\right)\nonumber\\
&\leq & N\Pr\left(d^2_1\II(|d_1|>x)>D\ln(N)\right)
\leq  N\Pr\left(d^2_1>\max\{D\ln(N), x^2\}\right)\nonumber\\
&\leq & N\Pr\left( \|\psi_{j_1k_1}\|_{\infty}|\varepsilon_{1, 0}|>1/2m_{11}m_{21}x\right)+N\Pr\left( \|\psi_{j_1k_1}\|_{\infty}|\xi_{1}|>1/2m_{11}m_{21}x\right)\nonumber\\
&\leq& 2N\exp\left\{-\frac{m_{11}^2m_{21}^2x^2}{8\|\psi_{j_1k_1}\|^2_{\infty}\EE[\varepsilon_{1,0}]^2}\right\}+N\exp\left\{-\frac{m_{11}^2m_{21}^2x^2}{8\|\psi_{j_1k_1}\|^2_{\infty}\EE[\xi_{1}]^2}\right\}\nonumber\\
& \leq & 2N\exp\left\{-\frac{m_{11}^2m_{21}^2\pi_o\ln(N)}{8\|\psi\|_{\infty}\EE[\varepsilon_{1,0}]^2}\right\}+2N\exp\left\{-\frac{m_{11}^2m_{21}^2\pi_o\ln(N)}{8\|\psi\|_{\infty}\EE[\xi_{1}]^2}\right\}.
\eeqn
Now, take $D=2[2\tau+1]\EE[\varepsilon_{1,0}]^2\frac{1}{m_{11}m_{21}}$ in \fr{p322}, $\pi_o=18[2\tau+1]\frac{\|\psi\|^2_{\infty}}{m_{11}m_{21}}\max\left\{\EE[\varepsilon_{1,0}]^2, \EE[\xi_{1}]^2\right\}$ in \fr{p321}, and $\gamma^2>\max\left\{\frac{288\tau}{m_{11}m_{21}}\|f\|^2_{\infty},\frac{288\tau}{m_{11}m_{21}}\|g\|^2_{\infty}, 256\tau^2\sigma^2\pi_o^2, 1152\tau D\sigma^2\right\}$ in  \fr{p1}, \fr{p2} and \fr{p31} completes the proof. \\
To prove \fr{Largdev-lm}, recall that for the case where $\sigma(t, x)$ is positive non-constant bivariate function, decomposition \fr{alef} will have four components, provided that $\int^1_0\int^1_0\psi_{j_1k_1}(t)\sigma(t, x)dtdx\neq 0$. The i.i.d parts and the martingale part will be dealt with in a similar manner as in the case when $\sigma(t, x)$ is constant. So, we will only prove large deviation result for the long-memory part. Indeed, since $\sum^N_{i=1}\varepsilon_{i; i-1}$ is a zero mean Gaussian random variable with variance \fr{longm},  by \fr{delta2} and \fr{thresh21}, we have
\beqn
\Pr\left(\left|\Delta_{22}\right|>1/12\lambda^{(1)}(j_1)\right)&\leq& 2 \exp\left\{-\frac{1/144\gamma_1^2\ln(N)}{2\int^1_0\int^1_0\frac{\psi_{j_1k_1}^2(t)\sigma^2(t, x)}{h_1(t)h_2(x)}dtdx\pi_{\alpha}}\right\}\nonumber\\
&\leq& 2 \exp\left\{-\frac{\gamma_1^2\ln(N)}{288\|\sigma\|^2_{\infty}m_{11}m_{21}\pi_{\alpha}}\right\}.
\eeqn
This completes the proof for a sufficiently large $\gamma_1$. The same way we can proof large deviation result for $| \widehat{\beta^{(2)}_{j_2k_2}}-\beta^{(2)}_{j_2k_2} |$. $ \Box$\\
{\bf Proof of Theorem \ref{th:upperbds-2}}.   First we will derive the upper bound for the $L^2$-risk for estimating the univariate function $f$, and the same procedure may be used for the estimation of $g$. Denote
\be \label{chijj}
\chi^{(1)}_{N}=\frac{\ln(N)}{A_1^2N}, \ \ 2^{J_{1o}}=[\chi_{N}]^{-\frac{1}{2s_1 + 1}},
\ee
 and note that with the choice of $J_1$ and $\lambda^{(1)}(j_1)$ given by \fr{J}  and \fr{thresh1}, respectively, the estimation error can be decomposed as $\mathbb{E} \| \widehat{f}_N-f \|^2
\leq\mathbb{E}_1 +\mathbb{E}_2 +\mathbb{E}_3+\mathbb{E}_4+\mathbb{E}_5+\mathbb{E}_6$, where
 \beqn  
\mathbb{E}_1&=&  \sum^{J_1-1}_{j_1=0}\sum^{2^{j_1}-1}_{k_1=0}\EE \left[\left| \widehat{ \beta^{(1)}_{j_1k_1}}-\beta^{(1)}_{j_1k_1}\right|^2 \II \left(   \left| \widehat{ \beta^{(1)}_{j_1k_1}}-\beta^{(1)}_{j_1k_1}  \right| > \frac{1}{2} \lambda^{(1)}(j_1)\right)  \right], \label{r21}\\
\mathbb{E}_2&=& \sum^{J_1-1}_{j_1=0} \sum^{2^{j_1}-1}_{k_1=0}\EE \left[\left| \widehat{ \beta^{(1)}_{j_1k_1}}-\beta^{(1)}_{j_1k_1}\right|^2 \II \left(  \left|  \beta^{(1)}_{j_1k_1} \right| >  \frac{1}{2}\lambda^{(1)}(j_1) \right)  \right],\label{r22}\\
\mathbb{E}_3&=& \sum^{J_1-1}_{j_1=0}\sum^{2^{j_1}-1}_{k_1=0}\left(\beta^{(1)}_{j_1k_1}\right)^2 \Pr \left( \left| \widehat{ \beta^{(1)}_{j_1k_1}}-\beta^{(1)}_{j_1k_1}  \right| > \frac{1}{2} \lambda^{(1)}(j_1) \right),\label{r31}\\
\mathbb{E}_4&=&\sum^{J_1-1}_{j_1=0}\sum^{2^{j_1}-1}_{k_1=0} \left(\beta^{(1)}_{j_1k_1}\right)^2\II \left(  \left|  \beta^{(1)}_{j_1k_1} \right| <  \frac{3}{2}\lambda_{l} \right).\label{r32}\\
\mathbb{E}_5&=& \sum^{\infty}_{j_1=J_1}\sum^{2^{j_1}-1}_{k_1=1} \left(\beta^{(1)}_{j_1k_1}\right)^2,\label{r33}\\
\mathbb{E}_6&=& \EE \left[\left| \widehat{\theta^{(1)}_{0, 0}} - \theta^{(1)}_{0, 0}\right|^2\right]. \label{r34}
\eeqn
Then, by  \fr{var-thet}, \fr{r34} becomes,
\beqn \label{e6}
\mathbb{E}_6=O\left(N^{-1}\right). 
\eeqn
Now, by \fr{eq11} and \fr{J}, \fr{r33} becomes 
\beqn \label{ee14}
\mathbb{E}_5
&=&O\left( \sum^{\infty}_{j_1=J_1} A_1^2 2^{-j_12s_1} \right)=O\left(A_1^2 \left[\chi_{N}\right]^{{2s_1}}\right)=O\left(A_1^2 \left[\chi_{N}\right]^{\frac{2s}{2s+1}}\right).
\eeqn
Now, combining $\mathbb{E}_1$ and $\mathbb{E}_3$, and applying Cauchy-Schwarz inequality, the moments property of the Gaussian, Lemma \ref{lem:Var} with the choice $\tau > 1$, \fr{eq11} and \fr{J}, yields 
\be
\mathbb{E}_1 +\mathbb{E}_3 = O  \left( \frac{1}{\ln(N)}\left[ \frac{1}{N}\right]^{{\tau}}+ A_1^2\left[\frac{\ln(N)}{N}\right]^{2s_1}\left[ \frac{1}{N}\right]^{{2\tau}} \right)=  o \left(  \frac{1}{N} \right). \label{r21r31}
\ee
Now, combining $\mathbb{E}_2$ and $\mathbb{E}_4$ and using condition \fr{eq11} yields
\be  \label{r22r32}
 \Delta=\mathbb{E}_2 +\mathbb{E}_4=  O \left( \sum^{J_1-1}_{j_1=0}\sum^{2^{j_1}-1}_{k_1=0}\min \left\{  \left( \beta^{(1)}_{j_1k_1}\right)^2,   A_1^2 \left[ \chi_{N}\right]\right\} \right).
\ee
Finally, $\Delta$ can be decomposed into the following components
\beqn  
 \Delta_1&=&  O \left( \sum^{{J_1}-1}_{j_1=j_{10}}A_1^2 2^{-2j_1s_1} \right)= O\left(A_1^2 \left[\chi_{N}\right]^{\frac{2s_1}{2s_1+1}}\right), \label{del1}\\
 \Delta_2&=&O \left( \sum^{j_{10}-1}_{j_1=0}A_1^2 2^{j_1} \left[ \chi_{N}\right]\right)=O \left( 2^{j_{10}}A_1^2  \left[ \chi_{N}\right]\right) =O\left(A_1^2 \left[\chi_{N}\right]^{\frac{2s_1}{2s_1+1}}\right). \label{del2}
 \label{del3}
\eeqn
Hence, combining \fr{e6}, \fr{ee14}, \fr{r21r31}, \fr{del1} and \fr{del2} yields
\beqn \label{sumofsq-f}
\mathbb{E} \| \widehat{f}_N-f \|^2=O\left(A_1^2 \left[\chi_{N}\right]^{\frac{2s_1}{2s_1+1}}\right).
\eeqn
The same way we can show that 
\beqn \label{sumofsq-g}
\mathbb{E} \| \widehat{g}_N-g \|^2=O\left(A_2^2 \left[\chi_{N}\right]^{\frac{2s_2}{2s_2+1}}\right).
\eeqn
Now, since wavelets functions are orthogonal to constants and scaling functions are not, the mean-squared error for estimating $U(t, x)$ by \fr{uest} is,  
\beqn \label{sumofsq}
\EE\|\widehat{U}(t, x)-U(t, x)\|_2^2&=&\EE\left|\widehat{\beta}_o-\beta_o\right|^2+\mathbb{E} \| \widehat{f}_N-f \|_2^2+\mathbb{E} \| \widehat{g}_N-g \|_2^2\nonumber\\
&+&2\|\psi^{(s)}_{0, 0}(t)\|_1\EE\left[(\widehat{\theta^{(1)}}_{0, 0}-{\theta^{(1)}}_{0, 0})(\widehat{\beta}_o-\beta_o)\right]+ 2\|\eta^{(s)}_{0, 0}(t)\|_1\EE\left[(\widehat{\theta^{(2)}}_{0, 0}-{\theta^{(2)}}_{0, 0})(\widehat{\beta}_o-\beta_o)\right].\nonumber\\
\eeqn
Hence, by Cauchy-Schwarz inequality, \fr{var-thet}, \fr{var-beto}, \fr{sumofsq-f} and \fr{sumofsq-g}, \fr{sumofsq} yields \fr{upperbds-2}. $\Box$\\

\section{Acknowledgements}

We thank Dr. German A. Schnaidt Grez and Dr. Brani Vidakovic for sharing their MATLAB codes on which we built our estimation algorithm.

\end{document}